\documentclass[10pt]{article}

\usepackage{amsmath}
\usepackage{amsmath, amsthm, amsfonts, amssymb, color}
\usepackage{mathrsfs}
\usepackage{latexsym,bm}
\usepackage{graphicx}
\usepackage[toc,page]{appendix}

\newtheorem{theorem}{Theorem}[section]%edit the theorems
\newtheorem{proposition}[theorem]{Proposition}

\newtheorem{lemma}[theorem]{Lemma}

\definecolor{wco}{rgb}{0.5,0.2,0.3}

\newcommand{\p}{\mathrm{P}}

\newcommand{\pp}{\mathbb{P}}
\newcommand{\qq}{\mathbb{Q}}

\newcommand{\mf}{\mathcal{M}(\mathbb{R})}

\newcommand{\R}{\mathbb{R}}
\newcommand{\e}{\mathrm{e}}
\newcommand{\C}{\mathbb{C}}
\newcommand{\re}{\mathrm{Re}}
\newcommand{\im}{\mathrm{Im}}
\numberwithin{equation}{section} %the number of the equation will be counted in one section
\begin{document}

\allowdisplaybreaks

\title{\bf Fluctuations of the additive martingales related to super-Brownian motion}
\author{
{\bf Ting Yang }\thanks{The research of this author is supported by NSFC (Grant Nos. 12271374 and 12371143).}
}

\date{}
\maketitle

\begin{abstract}
Let $(W_{t}(\lambda))_{t\ge 0}$, parametrized by $\lambda\in\R$, be the additive martingale related to a supercritical super-Brownian motion on the real line and let $W_{\infty}(\lambda)$ be its limit. Under a natural condition for the martingale limit to be non-degenerate, we investigate the rate at which the martingale approaches its limit. Indeed, assuming certain moment conditions on the branching mechanism, we show that the tail martingale $W_{\infty}(\lambda)-W_{t}(\lambda)$, properly normalized, converges in distribution to a non-degenerate random variable, and we identify the limit laws. We find that, for parameters with small absolute value, the fluctuations are affected by the behaviour of the branching mechanism $\psi$ around $0$. In fact, we prove that, in the case of small $|\lambda|$, when $\psi$ is secondly differentiable at $0$, the limit laws are scale mixtures of the standard normal laws, and when $\psi$ is `stable-like' near $0$ in some proper sense, the limit laws are scale mixtures of the stable laws. However, the effect of the branching mechanism is limited in the case of large $|\lambda|$. In the latter case, we show that the fluctuations and limit laws are determined by the limiting extremal process of the super-Brownian motion.
\end{abstract}

\medskip

\noindent\textbf{AMS 2020 Mathematics Subject Classification.} Primary 60J68; Secondary 60F05; 60G44

\medskip

\noindent\textbf{Keywords and Phrases.} Super-Brownian motion; Additive martingale; Rate of convergence; Central limit theorem.

\section{Introduction}

\subsection{Super-Brownian motion and the additive martingale}\label{sec1}

Throughout this paper, we use
``$:=$"
as a way of definition. The letters $c$ (with subscript) denote finite positive constants which may vary from place to place.
We use $\mathcal{B}_{b}(\R)$ (respectively,
$\mathcal{B}^{+}(\R)$) to denote
the space of bounded (respectively,
nonnegative) Borel functions on $\R$.
Suppose $\mf$ is the space of finite measures on $\mathbb{R}$ equipped with the topology of weak convergence.
For $\mu\in\mf$, we use the notation $\langle f,\mu\rangle:=\int_{\R}f(x)\mu(dx)$ and $\|\mu\|:=\langle 1,\mu\rangle$.
 We always assume $\{((B_{t})_{t\ge 0},\Pi_{x}):x\in\R\}$ is a standard Brownian motion on $\R$ and $P_{t}$ is the corresponding semigroup.

 Suppose $(X_{t})_{t\ge 0}$ is an $\mf$-valued Markov process. The distribution of $X$ is denoted by $\p_{\mu}$ if it is started
 at $\mu\in\mf$ at $t=0$.
 With abuse of notation, we also use $\p_{\mu}$ to denote the expectation with respect to $\p_{\mu}$. We write $\p$ as shorthand for $\p_{\delta_{0}}$.
 $X$ is called a
(supercritical) \textit{super-Brownian motion with branching mechanism $\psi$}
if for any
 $\mu\in \mf$, $f\in\mathcal{B}^{+}_{b}(\mathbb{R})$ and $t\ge 0$,
 \begin{equation}\nonumber
 \p_{\mu}\left[e^{-\langle f,X_{t}\rangle}\right]=e^{-\langle u_{f}(t, \cdot),\mu\rangle},
 \end{equation}
where
\begin{equation}\nonumber
u_{f}(t,x):=-\log \p_{\delta_{x}}\left(e^{-\langle f,X_{t}\rangle}\right)
\end{equation}
is
 the unique nonnegative locally bounded solution
 to the following integral equation:
 \begin{equation}\nonumber
 u_{f}(t,x)=P_{t}f(x)-\int_{0}^{t}P_{s}\left(\psi(u_{f}(t-s, \cdot))\right)(x)ds \quad
 \mbox{ for any }x\in \mathbb{R} \mbox{ and }t\ge 0.
 \end{equation}
Here $\psi$ is the so-called
branching mechanism,
 which takes the form
 \begin{equation}
 \psi(\lambda)=-a {\lambda}+b{\lambda}^{2}+\int_{(0,{\infty})}\left(e^{-{\lambda} y}-1+{\lambda} y\right)\pi(dy)\quad\mbox{ for }{\lambda}\ge 0,\label{2.1}
 \end{equation}
 with $a>0$, $b\ge 0$ and  $\pi(dy)$ being a measure concentrated on $(0,{\infty})$ such that
 $\int_{(0,+\infty)}\left(y\wedge y^{2}\right)\pi(dy)<+\infty $.
We note that $u_{f}(t,x)$ is also a solution to the partial differential equation
 \begin{equation}\label{eq1}
 \frac{\partial}{\partial t}u(t,x)=\frac{1}{2}\frac{\partial^{2}}{\partial x^{2}} u(t,x)-\psi(u(t,x))
 \end{equation}
 with initial condition $u(0,x)=f(x)$.
\eqref{eq1} is referred to as the Kolmogorov-Petrovsky-Piscounov (K-P-P) equation.
We refer to \cite{Li} for the existence of super-Brownian motions.

We assume $\psi(+\infty)=+\infty$ which assures that the event $\{\lim_{t\to+\infty}\|X_{t}\|=0\}$ will occur with positive probability. Indeed, there exists some $\lambda^{*}\in (0,+\infty)$ satisfying $\psi(\lambda^{*})=0$ such that
$$\p_{\mu}\left(\lim_{t\to+\infty}\|X_{t}\|=0\right)=\e^{-\lambda^{*}\|\mu\|}\quad\forall \mu\in\mf.$$
For $\lambda\in\R$, define $e_{\lambda}(x):=\mathrm{e}^{-\lambda x}$ for $x\in\R$, and
\begin{equation*}
W_{t}(\lambda):=\mathrm{e}^{-\left(\frac{\lambda^{2}}{2}+a\right)t}\langle e_{\lambda},X_{t}\rangle\quad\forall t\ge 0.
\end{equation*}
It is well-known that for every $\lambda\in\R$, $(W_{t}(\lambda))_{t\ge 0}$ forms a nonnegative $\p$-martingale with respect to the natural filtration $\sigma(X_{s}:s\le t)$, and thus converges a.s. to a random variable $W_{\infty}(\lambda)$. This martingale is called the additive martingale of the super-Brownian motion. Interest of this martingale is stimulated by its close connection with the travelling wave solution to the K-P-P equation, see, for example, \cite{KLMR} and the references therein.
The counterpart of this martingale in the setting of branching random walk (BRW) is so called the Biggins' martingale, which plays a key role in the study of BRW. The convergence of Biggins' martingale were studied in \cite{B1,B2,ILL,KM}.

Assuming mild conditions on the branching mechanism, \cite{KLMR} established the necessary and sufficient conditions for the additive martingale to converge to a non-degenerate limit.
By \cite[Theorem 2.4]{KLMR}, assume that
\begin{description}
\item{(A0)} $$\int^{+\infty}\frac{1}{\sqrt{\int_{\lambda^{*}}^{\xi}\psi(u)du}}d\xi<+\infty,$$
\end{description}
then $W_{\infty}(\lambda)$ is non-degenerate if and only if
\begin{equation}\label{condi2}
|\lambda|<\sqrt{2a}\mbox{ and }\int_{(1,+\infty)}r\log r\pi(dr)<+\infty.
\end{equation}
In this case $W_{t}(\lambda)\to W_{\infty}(\lambda)$ in $L^{1}(\p)$ and $W_{\infty}(\lambda)>0$ $\p$-a.s. on the survival event $\mathcal{E}^{c}:=\{\|X_{t}\|>0\mbox{ for all }t>0\}$.
For the boundary case $|\lambda|=\sqrt{2a}$, $W_{t}(\lambda)$ converges a.s. to $0$ and a natural object to study is the so-called derivative martingale
$$\partial W^{\pm}_{t}:=\e^{-2a t}\langle(\sqrt{2a}t\pm\cdot)e_{\pm \sqrt{2a}},X_{t} \rangle\quad\forall t\ge 0.$$
Due to \cite[Theorem 2.4]{KLMR}, assume (A0) and that
\begin{equation}\label{condi3}
\int_{(1,+\infty)}r(\log r)^{2}\pi(dr)<+\infty,
\end{equation}
$\partial W^{\pm}_{t}\to \partial W_{\infty}^{\pm}$ $\p$-a.s. as $t\to+\infty$ and
$\partial W_{\infty}^{\pm}>0$ a.s. on $\mathcal{E}^{c}$.
The limits of the derivative martingale appear in many limit theorems of super-Brownian motions. For example,
it is proved by \cite[Theorem 1.1]{HRS} that under (A0) and \eqref{condi3},
\begin{equation}\label{HRSthem1.1}
\sqrt{t}W_{t}(\sqrt{2a})\to \sqrt{\frac{2}{\pi}}\partial W_{\infty}^{+}\mbox{ in probability as }t\to+\infty.
\end{equation}
Due to symmetry, the case $\lambda=-\sqrt{2a}$ can be analyzed similarly and one has $\sqrt{t}W_{t}(-\sqrt{2a})\to \sqrt{\frac{2}{\pi}}\partial W_{\infty}^{-}$ in probability as $t\to+\infty$.
By \cite{RSZ,RYZ}, $\partial W^{-}_{\infty}$ also appears in the weak limit of the extremal process of the super-Brownian motion (see, Section \ref{sec2.1} below, for a short discussion).

It is natural to ask how fast the martingale $W_{t}(\lambda)$ approaches its limit $W_{\infty}(\lambda)$ as $t\to+\infty$.
For the boundary case $|\lambda|=\sqrt{2a}$, $W_{\infty}(\lambda)$ vanishes and \cite[Theorem 1.1]{HRS} can be viewed as a statement on the rate of convergence.
One is thus interested in the convergence rate of the tail martingale $W_{\infty}(\lambda)-W_{t}(\lambda)$ when $W_{\infty}(\lambda)$ is non-degenerate.
There are many earlier results on the rate of convergence for the tail martingales related to various branching models, including Galton-Watson process, BRW, random trees and etc. We refer to \cite[Page 1182]{IK} and \cite[Page 2454]{IKM2} for a detailed account of the literature.

In the setting of BRW, the rate of convergence for Biggins' martingale has been obtained in various sense.
See, \cite{I,IM} for the rate of almost sure convergence, \cite{AIPR} for the rate of $L^{p}$-convergence, and \cite{IK,IKM1,IKM2} for the rate in the form of a central limit theorem.
In fact, \cite{IK,IKM1} obtained the rate in the form of a CLT for real-valued and sufficiently small parameters, and then \cite{IKM2} extended the results in \cite{IK} to all complex-valued parameters.

For supercritical superprocesses, the linear statistics and related convergence rates are objects of active studies. We refer to \cite{RSSZ,RSZ0} (and the references therein) for the rates in the form of CLTs, and refer to \cite{LRS} for the rates of almost sure convergence and $L^{p}$-convergence. Although the additive martingale is a special case of linear statistic, the conditions imposed in the aforementioned papers exclude the case of spatial-homogeneous super-Brownian motions.

We also notice that a recent paper \cite{HRS1} considered the rate of almost sure convergence for tail martingales related to branching Brownian motion (BBM).
Another related work can be found in recent papers \cite{MP} and \cite{HRS2}, where fluctuations of the derivative martingale were studied, respectively, for BBM and BRW.

In this paper we shall investigate the rate of convergence in the form of a CLT for the tail martingale related to super-Brownian motion. Our two main results, Theorem \ref{them3.1} and Theorem \ref{them5.1}, show that, assuming a natural sufficient condition for the martingale limit to be non-degenerate and some additional moment conditions on the branching mechanism, the tail martingale $W_{\infty}(\lambda)-W_{t}(\lambda)$, properly normalized, converges in distribution to a non-degenerate random variable, and we identify the limit laws.

\subsection{Stable distribution}
Before we state the main results, we introduce some known facts on the stable distributions. A standard reference for the facts given in this subsection is \cite{Kyprianou2}.

For a complex valued function $f:\R\to \C$ such that $f(0)=1$ and $f(x)\not=0$ for $x\in\R$, there exists a unique continuous function $g$ from $\R$ to $\C$ such that $g(0)=0$ and $f(x)=\e^{g(x)}$ for all $x\in\R$. This function $g$ is denoted by $\log f$ and is called the distinguished logarithm of $f$.
An $\R$-valued random variable $\xi$ is said to have a stable distribution if for any $a_{1},a_{2}>0$ there are $c>0$ and $b\in \R$ such that
$$a_{1}\xi_{1}+a_{2}\xi_{2}\stackrel{d}{=}c\xi+b$$
where $\xi_{1},\xi_{2}$ are independent copies of $\xi$.
It is known that the distinguished logarithm of the characteristic function of a stable random variable $\xi$  has the following form
\begin{equation*}
\log\mathrm{E}[\e^{i\theta \xi}]=
\begin{cases}
-\nu^{\alpha}|\theta|^{\alpha}(1-i\eta \mathrm{sign}(\theta)\tan(\frac{\pi}{2}\alpha))+i\mu\theta,&\quad \alpha\not=1;\\
-\nu|\theta|(1+i\eta\frac{2}{\pi}\mathrm{sign}(\theta)\log|\theta|)+i\mu\theta,& \quad \alpha=1,
\end{cases}
\end{equation*}
for some constants $\alpha\in (0,2)$, $\eta\in [-1,1]$, $\nu>0$ and $\mu\in\R$.  In particular, $\xi$ is said to have a strictly stable distribution if
\begin{equation*}
\log\mathrm{E}[\e^{i\theta \xi}]=
\begin{cases}
-\nu^{\alpha}|\theta|^{\alpha}(1-i\eta \mathrm{sign}(\theta)\tan(\frac{\pi}{2}\alpha)),&\quad \alpha\not=1;\\
-\nu|\theta|+i\mu\theta,& \quad \alpha=1,
\end{cases}
\end{equation*}
where the parameters ranges for $\alpha,\nu,\eta,\mu$ are as above. In this case, one has for all $n\ge 1$,
$$\xi_{1}+\xi_{2}+\cdots+\xi_{n}\stackrel{d}{=}n^{1/\alpha}\xi$$
where $\xi_{1},\cdots,\xi_{n}$ are independent copies of $\xi$.
We use $\mathcal{S}_{\alpha}(\nu,\eta,\mu)$ to denote the stable distribution (also called the $\alpha$-stable law).
For example, if $\xi$ has the characteristic function $\mathrm{E}[\e^{i\theta \xi}]=\e^{(-i\theta)^{1+\beta}}$ for some $\beta\in (0,1)$, then $\xi$ has the strictly $(1+\beta)$-stable law $\mathcal{S}_{1+\beta}(\sin(\frac{\pi}{2}\beta)^{1/(1+\beta)},1,0)$.

\subsection{The main result}
We observe that up to a space-time scaling transform, the branching mechanism $\psi$ with $\psi(+\infty)=+\infty$ can always be assumed to satisfy that
\begin{equation}\label{condi0}
\psi'(0+)=-1\mbox{ and }\lambda^{*}=1.
\end{equation}
In fact, if $u(t,x)$ is a solution to \eqref{eq1} then $u(a^{-1}t,a^{-1/2}x)/\lambda^{*}$ is a solution to \eqref{eq1} with $\psi$ replaced by $\widetilde{\psi}$ where $\widetilde{\psi}(\lambda)=\psi(\lambda^{*}\lambda)/a\lambda^{*}$ satisfies that $\widetilde{\psi}'(0+)=-1$ and $\tilde{\psi}(1)=0$. This implies that
for a super-Brownian motion
$(X_{t})_{t\ge 0}$ with branching mechanism $\psi$, if we define the random measures $\tilde{X}_{t}$ by
\begin{equation*}
\langle f,\widetilde{X}_{t}\rangle =\lambda^{*}\langle f(a^{1/2}\cdot),X_{a^{-1}t}\rangle \quad\forall t\ge 0,\ f\in\mathcal{B}_{b}^{+}(\mathbb{R}),
\end{equation*}
then $(\widetilde{X}_{t})_{t\ge 0}$ is a
super-Brownian motion with branching mechanism $\widetilde{\psi}$.
Let $\widetilde{W}_{t}(\lambda)$ denote the additive martingale with respect to $(\widetilde{X}_{t})_{t\ge 0}$. We necessarily have $\widetilde{W}_{t}(\lambda)=\lambda^{*}W_{a^{-1}t}(a^{1/2}\lambda)$ for all $\lambda\in\R$ and $t\ge 0$.
Therefore, in the rest of this paper the branching mechanism $\psi$ is assumed to satisfy \eqref{condi0}, which will simplify computations and notations later on.

Let $\bar{\pi}(r):=\int_{(r,+\infty)}\pi(dy)$ for $r>0$  and
$$\psi_{0}(\lambda):=\lambda+\psi(\lambda)=b\lambda^{2}+\int_{(0,+\infty)}\left(\e^{-\lambda y}-1+\lambda y\right)\pi(dy)\quad\forall \lambda\ge 0.$$
We will assume the branching mechanism $\psi$ satisfies (A0) and either one of the following conditions.
\begin{description}
\item{(A1)} $\int_{1}^{+\infty}r\bar{\pi}(r)dr<+\infty$.
\item{(A2)} There are constants $0<\beta<\beta+\delta<1$, $A_{1}>0$ and $c_{1}\ge 0$ such that for $r>0$ sufficiently large
$$\left|\bar{\pi}(r)-\frac{A_{1}}{r^{1+\beta}}\right|\le \frac{c_{1}}{r^{1+\beta+\delta}}.$$
\end{description}

We remark that (A1) is equivalent to requiring that
$$\int_{(1,+\infty)}r^{2}\pi(dr)<+\infty.$$
This is because by Fubini's theorem $\int_{1}^{+\infty}r\bar{\pi}(r)dr=\int_{(1,+\infty)}\frac{1}{2}(y^{2}-1)\pi(dy)$. Indeed, (A1) requires that $\psi$ is secondly differentiable at $0$, and thus assures the finiteness of variance of the super-Brownian motion.
On the other hand, (A2) requires that $\psi_{0}$ is not far away from $\lambda^{1+\beta}$ near $0$. In fact, by Lemma \ref{lem5.1} below, one can verify that under (A2) $\psi_{0}(\lambda)/\lambda^{1+\beta}\to c$ as $\lambda\to 0+$ for some constant $c>0$. Hence (A2) allows infinite variance of the process.

We shall also assume on occasion that
\begin{description}
\item{(A3)} There are constants $\gamma\in (0,1]$ and $c_{2}>0$ such that for $\lambda>0$ sufficiently large,
$$\psi_{0}(\lambda)\ge c_{2}\lambda^{1+\gamma}.$$
\end{description}
Obviously (A3) implies (A0). (A3) is automatically true if $b>0$. When $b=0$, this condition is equivalent to requiring that $\int_{0}^{1}\bar{\pi}(r)(r\wedge \lambda)dr\ge c_{3}\lambda^{1-\gamma}$ for $\lambda>0$ sufficiently small (see, Lemma \ref{lemA.0} below).

For example, if we take $\psi(\lambda)=-\lambda+(1+\beta)^{-1}\lambda^{1+\beta}$ for some $\beta\in (0,1]$, then (A1)(A3) are satisfied if $\beta=1$,  and (A2)(A3) are satisfied if $\beta\in (0,1)$.

The results are stated separately for branching mechanisms satisfying (A1) and (A2).
\begin{theorem}\label{them3.1}
Suppose (A0) and (A1) hold.
\begin{description}
\item{(i)} For $|\lambda|<\frac{\sqrt{2}}{2}$, under $\p$,
$$\e^{\frac{1-\lambda^{2}}{2}t}(W_{\infty}(\lambda)-W_{t}(\lambda))\stackrel{d}{\to}\sqrt{\frac{\sigma^{2}}{1-\lambda^{2}}}\cdot\sqrt{W_{\infty}(2\lambda)}\cdot U_{1}\mbox{ as }t\to+\infty,$$
where $\sigma^{2}:=2b+\int_{(0,+\infty)}r^{2}\pi(dr)$, $W_{\infty}(2\lambda)$ is non-degenerate, and $U_{1}$ is a standard normal random variable which is independent of $W_{\infty}(2\lambda)$.
\end{description}
Assume, in addition, (A3) holds. Then the following is true.
\begin{description}
\item{(ii)} For $\lambda=\pm\frac{\sqrt{2}}{2}$, under $\p$,
$$t^{1/4}\e^{t/4}(W_{\infty}(\lambda)-W_{t}(\lambda))\stackrel{d}{\to}2^{3/4}\pi^{-1/4}\sigma\cdot\sqrt{\partial W_{\infty}^{\pm}}\cdot U_{2}\mbox{ as }t\to+\infty,$$
where $U_{2}$ is a standard normal random variable independent of $\partial W_{\infty}^{\pm}$.
\item{(iii)} For $\frac{\sqrt{2}}{2}<|\lambda|<\sqrt{2}$, there exists a non-degenerate random variable $\eta_{\lambda}$ such that
$$t^{\frac{3|\lambda|}{2\sqrt{2}}}\e^{\frac{1}{2}\left(\sqrt{2}-|\lambda|\right)^{2}t}\left(W_{\infty}(\lambda)-W_{t}(\lambda)\right)\stackrel{d}{\to} \eta_{\lambda}\mbox{ as }t\to+\infty.$$
Moreover, the characteristic function of $\eta_{\lambda}$ is given in \eqref{eq:cf}.
\end{description}
\end{theorem}

\begin{theorem}\label{them5.1}
Assume (A0) and (A2) hold.
\begin{description}
\item{(i)} For $|\lambda|<\frac{\sqrt{2}}{1+\beta}$, under $\p$,
$$\e^{\frac{\beta}{2}\left(\frac{2}{1+\beta}-\lambda^{2}\right)t}(W_{\infty}(\lambda)-W_{t}(\lambda))\stackrel{d}{\to}
c_{\lambda,\beta}^{\frac{1}{1+\beta}}\cdot W_{\infty}((1+\beta)\lambda)^{\frac{1}{1+\beta}}\cdot V_{1},$$
Where $c_{\lambda,\beta}:=A_{1}\beta^{-2}\Gamma(1-\beta)(1-\frac{1}{2}\lambda^{2}(1+\beta))^{-1}$, $W_{\infty}((1+\beta)\lambda)$ is non-degenerate, and $V_{1}$ is a strictly $(1+\beta)$-stable random variable which is independent of $W_{\infty}((1+\beta)\lambda))$.
\end{description}
Assume, in addition, (A3) holds. Then the following is true.
\begin{description}
\item{(ii)} For $\lambda=\pm\frac{\sqrt{2}}{1+\beta}$, under $\p$,
$$t^{\frac{1}{2(1+\beta)}}\e^{\frac{\beta^{2}}{(1+\beta)^{2}}t}(W_{\infty}(\lambda)-W_{t}(\lambda))\stackrel{d}{\to}c_{\beta}^{\frac{1}{1+\beta}}\cdot (\partial W^{\pm}_{\infty})^{\frac{1}{1+\beta}}\cdot V_{2}\mbox{ as }t\to+\infty,$$
Where $c_{\beta}:=A_{1}\beta^{-3}(1+\beta)\Gamma(1-\beta)\sqrt{2/\pi}$, and $V_{2}$ is a strictly $(1+\beta)$-stable random variable which is independent of $\partial W_{\infty}^{\pm}$.
\item{(iii)} For $\frac{\sqrt{2}}{1+\beta}<|\lambda|<\sqrt{2}$, there exists a non-degenerate random variable $\eta_{\lambda}$ such that
$$t^{\frac{3|\lambda|}{2\sqrt{2}}}\e^{\frac{1}{2}\left(\sqrt{2}-|\lambda|\right)^{2}t}\left(W_{\infty}(\lambda)-W_{t}(\lambda)\right)\stackrel{d}{\to} \eta_{\lambda}\mbox{ as }t\to+\infty.$$
Moreover, the characteristic function of $\eta_{\lambda}$ is given in \eqref{eq:cf}.
\end{description}
\end{theorem}

Condition (A1) guarantees that $W_{\infty}(\lambda)$ is the $L^{2}(\p)$-limit of $W_{t}(\lambda)$ for all $|\lambda|\le \frac{\sqrt{2}}{2}$ (cf. \cite[Theorem 1.1]{KyM}).
So our Theorem \ref{them3.1}(i) and (ii) can be viewed as the super-Brownian motion counterparts of \cite[Theorem 2.2 and Theorem 2.3]{IKM2} in the case of real-valued and sufficient small parameters. In the setting of complex BBM random energy model, which includes BBM as a special case, an analog of Theorem \ref{them3.1}(i) was obtained in \cite{HK}.

While the fluctuations in Theorem \ref{them3.1}(i)and (ii) are Gaussian, Theorem \ref{them5.1}(i) and (ii) provide sufficient conditions for the tail martingale to exhibit stable-like fluctuations.
Analogues of Theorem \ref{them5.1}(i) are obtained in \cite{Heyde} and \cite{IKM2}, respectively, for the Galton-Watson process and BRW.

We remark that, in the case of large parameters, the weak convergence results appear in Theorem \ref{them3.1}(iii) and Theorem \ref{them5.1}(iii) hold
under weaker conditions. Indeed, our Proposition \ref{prop2} below shows that the convergence result holds for all $\frac{\sqrt{2}}{2}<|\lambda|<\sqrt{2}$, under (A3) and the assumption that $\int_{(1,+\infty)}r^{p}\pi(dr)<+\infty$ for some $p\in (\sqrt{2}/|\lambda|,2/\lambda^{2})$ and $p\le 2$. Neither (A1) nor (A2) is required there.
We also remark that the law of $\eta_{\lambda}$, whose characteristic function is given in \eqref{eq:cf} below, is determined by the limiting extremal process of super-Brownian motion.
Proposition \ref{prop2} can be viewed as the counterpart of \cite[Theorem 2.5]{IKM2} in the case of real-valued and large parameters.

\section{Preliminary of proofs}
In this section we gather all the auxiliary facts that will be used in the proofs of main results.

\subsection{Extremal process of super-Brownian motion}\label{sec2.1}
For $x\in\R$ and a measure $\mu$ on $\R$, we use $\mu+x$ to denote the measure induced by the shift operator $f\mapsto f(\cdot+x)$ for $f\in\mathcal{B}^{+}(\R)$. That is to say, $\langle f,\mu+x\rangle=\langle f(\cdot+x),\mu\rangle$ for all $f\in\mathcal{B}^{+}(\R)$.
Let $\max\mu$ denote the supremum of the support of $\mu$, i.e., $\max\mu:=\sup\{x\in \R:\ \langle 1_{(-\infty,x]},\mu\rangle>0\}$.
It is proved in \cite{RYZ} that, assume the branching mechanism $\psi$ satisfies \eqref{condi0} and that
\begin{equation}\label{condi4}
\int_{(1,+\infty)}r^{1+\alpha}\pi(dr)<+\infty\mbox{ for some }\alpha\in (0,1),
\end{equation}
then under $\p$, the extremal process of $X_{t}$, defined by
$$\mathcal{E}_{t}:=X_{t}-m(t),\mbox{ where }m(t):=\sqrt{2}t-\frac{3}{2\sqrt{2}}\log t$$
converges in distribution to a random Radon measure $\mathcal{E}_{\infty}$ on $\R$.
We point out that the above convergence result is also obtained in \cite{RSZ} assuming a slightly different moment condition on the branching mechanism.

In particular, if (A3) holds in addition, then $\max\mathcal{E}_{t}\stackrel{d}{\to}\max\mathcal{E}_{\infty}$ as $t\to+\infty$ and $\max \mathcal{E}_{\infty}$ has a randomly shifted Gumbel distribution.
In this case, the limiting extremal process $\mathcal{E}_{\infty}$ can be constructed as a decorated Poisson random measure, where the Poisson random measure has an exponential intensity, and each atom is decorated by an independent copy of an auxiliary random measure.
To be more specific, $\mathcal{E}_{\infty}$ can be constructed as follows: There exist a positive constant $\tilde{c}_{0}$ and a random Radon measure $\tilde{\triangle}$ supported on $(-\infty,0]$ such that,
given $\partial W_{\infty}^{-}$, let $\{e_{i}:i\ge 1\}$ be the atoms of a Poisson point process on $\R$ with intensity $\tilde{c}_{0}\partial W_{\infty}^{-}\sqrt{2}\e^{-\sqrt{2}y}dy$ and $\{\tilde{\triangle}_{i}:i\ge 1\}$ be an independent sequence of independent copies of $\tilde{\triangle}$, then
$$\mathcal{E}_{\infty}\stackrel{d}{=}\sum_{i\ge 1}(\tilde{\triangle}_{i}+e_{i}).$$
Let $\mathcal{H}:=\{f\in\mathcal{B}^{+}_{b}(\R):\ \int_{0}^{+\infty}y\e^{\sqrt{2}y}f(-y)dy<+\infty\}$.
The next lemma refines \cite[Theorem 2.10]{RYZ} where the result is obtained for $f\in\mathcal{H}$.

\begin{lemma}\label{lem4.1}
Assume (A3) and \eqref{condi4} hold.
If $f:\R\to \R^{+}$ is a nonnegative locally bounded function such that $f1_{(-\infty,K]}\in\mathcal{H}$ for some $K\in \R$, then $\langle f,\mathcal{E}_{\infty}\rangle$ is a.s. finite, and $\langle f,\mathcal{E}_{t}\rangle\stackrel{d}{\rightarrow}\langle f,\mathcal{E}_{\infty}\rangle$ as $t\to+\infty$.
\end{lemma}

\proof We write $f=f_{1}+f_{2}$ where $f_{1}(x)=f(x)1_{(-\infty,K]}(x)$. Since $f_{1}\in\mathcal{H}$, $\langle f_{1},\mathcal{E}_{\infty}\rangle$ is a.s. finite by \cite[Theorem 2.10]{RYZ}. To prove $\langle f,\mathcal{E}_{\infty}\rangle$ is a.s. finite, we only need to show that $\langle f_{2},\mathcal{E}_{\infty}\rangle$ is a.s. finite. By the construction of $\mathcal{E}_{\infty}$, we have
\begin{equation}\nonumber
\langle f_{2},\mathcal{E}_{\infty}\rangle\stackrel{d}{=}\sum_{i\ge 1}\langle f_{2}(\cdot+e_{i}),\tilde{\triangle}_{i}\rangle.
\end{equation}
Recall that $\tilde{\triangle}_{i}$ is supported on $(-\infty,0]$ and $f_{2}$ is supported on $[K,+\infty)$. So for $e_{i}<K$, $\langle f_{2}(\cdot+e_{i}),\tilde{\triangle}_{i}\rangle=0$. Thus
\begin{equation}
\langle f_{2},\mathcal{E}_{\infty}\rangle\stackrel{d}{=}\sum_{e_{i}\ge K}\int_{(K-e_{i},0]}f(y+e_{i})\tilde{\triangle}_{i}(dy).\label{lem4.1.1}
\end{equation}
Given $\partial W_{\infty}^{-}$, $\sum_{e_{i}\ge K}\delta_{e_{i}}$ is a Poisson random variable with the parameter $\tilde{c}_{0}\partial W_{\infty}^{-}\int_{K}^{+\infty}\sqrt{2}\e^{-\sqrt{2}y}dy<+\infty$. Thus the right hand side of \eqref{lem4.1.1} is a sum of finite terms.
Hence we prove $\langle f_{2},\mathcal{E}_{\infty}\rangle$ is a.s. finite.

For $n\in \mathbb{N}$, let $g_{n}(x):=f(x)1_{(-\infty,n]}(x)$. Obviously $\langle g_{n},\mathcal{E}_{\infty}\rangle \uparrow \langle f,\mathcal{E}_{\infty}\rangle$ a.s. as $n\to+\infty$. By our assumptions, $g_{n}\in\mathcal{H}$ for every $n\in\mathbb{N}$. Thus by \cite[Theorem 2.10]{RYZ} $\langle g_{n},\mathcal{E}_{t}\rangle \stackrel{d}{\rightarrow}\langle g_{n},\mathcal{E}_{\infty}\rangle$ as $t\to+\infty$ for each $n\in\mathbb{N}$. By \cite[Theorem 3.2]{Billingsley} the convergence $\langle f,\mathcal{E}_{t}\rangle\stackrel{d}{\to}\langle f,\mathcal{E}_{\infty}\rangle$ holds once we prove that for any $\epsilon>0$,
\begin{equation}\label{lem4.1.2}
\lim_{n\to+\infty}\limsup_{t\to+\infty}\pp\left(\left|\langle g_{n},\mathcal{E}_{t}\rangle-\langle f,\mathcal{E}_{t}\rangle\right|>\epsilon\right)=0.
\end{equation}
In fact, we have
\begin{equation}\nonumber
\pp\left(\left|\langle g_{n}-f,\mathcal{E}_{t}\rangle\right|>\epsilon\right)=\pp\left(\langle f1_{(n,+\infty)},\mathcal{E}_{t}\rangle>\epsilon\right)\le\pp(\max \mathcal{E}_{t}>n).
\end{equation}
Recall that under our assumptions, $\max \mathcal{E}_{t}$ converges in distribution to some random variable. Thus $\lim_{n\to+\infty}\lim_{t\to+\infty}\pp\left(\max \mathcal{E}_{t}>n\right)=0$, and \eqref{lem4.1.2} follows.\qed

\subsection{Skeleton BBM and related facts}

 Following \cite{BKM}, we give the so-called skeleton construction of a (supercritical) super-Brownian motion,
 which provides a pathwise description of
a super-Brownian motion in terms of immigrations along a BBM called the skeleton.

\begin{proposition}\label{prop1}
For every $\mu\in \mf$ and
every finite point measure $\nu$ on $\R$,
there exists a probability space with probability measure $\pp_{\mu,\nu}$ that carries three processes: $(Z_{t})_{t\ge 0}$, $(I_{t})_{t\ge 0}$ and $(X^{*}_{t})_{t\ge0}$, where
$\left((Z_{t})_{t\ge 0};\pp_{\mu,\nu}\right)$ is a BBM with  branching rate $q>0$ and offspring distribution $\{p_{k}:k\ge 2\}$, which are
uniquely determined by
\begin{equation}\label{eq2}
q(F(s)-s)=\psi(1-s)\quad\forall s\in [0,1],
\end{equation}
here $F(s):=\sum_{k=2}^{+\infty}p_{k}s^{k}$, and $\pp_{\mu,\nu}(Z_{0}=\nu)=1$;
$\left((X^{*}_{t})_{t\ge 0};\pp_{\mu,\nu}\right)$ is a (subcritical)
super-Brownian motion
with branching mechanism
$\psi^{*}(\lambda):=\psi(\lambda+1)$ and
$\pp_{\mu,\nu}(X^{*}_{0}=\mu)=1$; $((I_{t})_{t\ge 0},\pp_{\mu,\nu})$ is an $\mf$-valued process with $\pp_{\mu,\nu}\left(I_{0}=0\right)=1$, which denotes the immigration at time $t$ that occurred along the skeleton $Z$. Under $\pp_{\mu,\nu}$, both $Z$ and $I$ are independent of $X^{*}$.

If $\pp_{\mu}$ denotes the measure $\pp_{\mu,\nu}$ with $\nu$ replaced by a Poisson random measure with intensity $\mu(dx)$,
then $\left(
\widehat X :=
X^{*}+I;\pp_{\mu}\right)$ is Markovian and
has the same distribution as
$(X;\p_{\mu})$.
Moreover, under $\pp_{\mu}$, given $\widehat  X_{t}$,
the measure $Z_{t}$ is a Poisson random measure with intensity
$\widehat X_{t}(dx)$.
\end{proposition}

We shall only need the above defining properties in order to proceed the martingale fluctuations. We refer the readers to \cite{BKM} for further detailed information on the distributional properties of the skeleton space.

Since
$({\widehat X};\pp_{\mu})$
is equal in distribution to the super-Brownian motion $(X;\p_{\mu})$, we may work on this skeleton space whenever it is convenient.
 For notational simplification, we write $\pp$ as short hand for $\pp_{\delta_{0}}$, and we will abuse the notation and denote $\widehat X$ by $X$.
We will refer to $(Z_{t})_{t\ge 0}$ as the skeleton BBM of $X$.
Since the distributions of $Z$ under $\pp_{\mu,\nu}$ do not depend on $\mu$, we write $\qq_{\nu}$ for $\pp_{\mu,\nu}$ to simplify the notation.

We define the semigroup $P^{c}_{t}$ for $c\in\R$ by
$$P^{c}_{t}f(x):=\e^{c t}P_{t}f(x)\quad\forall f\in\mathcal{B}^{+}(\R).$$
For every $x\in\R$, the
first moments of $X_{t}$ and $Z_{t}$
can be expressed as follows:  For $x\in\R$, $f\in\mathcal{B}_{b}^{+}(\R)$ and $t\ge 0$,
\begin{equation}\nonumber
\pp_{\delta_{x}}\left[\langle f,X_{t}\rangle\right]=\qq_{\delta_{x}}\left[\langle f,Z_{t}\rangle\right]=P^{1}_{t}f(x).
\end{equation}
Using the above equality, one can easily verify that for $\lambda\in\R$ and $c(\lambda):=\frac{1}{2}\lambda^{2}+1$,
$$Z_{t}(\lambda):=\e^{-c(\lambda)t}\langle e_{\lambda},Z_{t}\rangle\quad\forall t\ge 0$$
is a nonnegative $\qq_{\delta_{x}}$-martingale for every $x\in\R$, and thus converges a.s. to a finite random variable $Z_{\infty}(\lambda)$.
$Z_{t}(\lambda)$ is referred to as the additive martingale of the skeleton BBM.

Let $L$ be an integer-valued random variable with the generating function $F(s)$ given by \eqref{eq2}. By \cite[Theorem 1]{Kyprianou}, the martingale limit $Z_{\infty}(\lambda)$ is non-zero with positive probability if and only if $|\lambda|<\sqrt{2}$ and $\mathrm{E}[L\log^{+}L]<+\infty$. When the latter holds, $Z_{t}(\lambda)\to Z_{\infty}(\lambda)$ in $L^{1}(\qq_{\delta_{x}})$ for every $x\in\R$. The following $L^{p}$-convergence result is from \cite[Theorem 1.3]{HH}.
\begin{lemma}\label{lem2.0}
Suppose $|\lambda|<\sqrt{2}$. For $p\in (1,\frac{2}{\lambda^{2}})$ and $p\le 2$ such that $\mathrm{E}[L^{p}]<+\infty$,  $Z_{t}(\lambda)$ converges to $Z_{\infty}(\lambda)$
$\qq_{\delta_{x}}$-a.s. and in $L^{p}(\qq_{\delta_{x}})$ for every $x\in\R$.
\end{lemma}

\begin{lemma}\label{lem2.1}
Suppose (A0) and \eqref{condi2} hold. For $|\lambda|<\sqrt{2}$,
$$W_{\infty}(\lambda)=Z_{\infty}(\lambda)\quad\pp\mbox{-a.s.}$$
\end{lemma}
\proof By \cite[Theorem 2.4]{KLMR}, under our assumptions, $W_{\infty}(\lambda)$ is the $L^{1}(\pp)$-limit of $W_{t}(\lambda)$. Let $\mathcal{F}_{t}$ be the $\sigma$-field generated by $Z$, $X^{*}$ and $I$ up to time $t$. Let $s>0$. Then
\begin{eqnarray}
\pp\left[|\pp\left[W_{t+s}(\lambda)|\mathcal{F}_{t}\right]-W_{\infty}(\lambda)|\right]
&\le&\pp\left[|\pp\left[W_{t+s}(\lambda)|\mathcal{F}_{t}\right]-\pp\left[W_{\infty}(\lambda)|\mathcal{F}_{t}\right]|\right]
+\pp\left[|\pp\left[W_{\infty}(\lambda)|\mathcal{F}_{t}\right]-\pp\left[W_{t}(\lambda)|\mathcal{F}_{t}\right]|\right]\nonumber\\
&&\quad+\pp[|W_{t}(\lambda)-W_{\infty}(\lambda)|]\nonumber\\
&\le&\pp\left[|W_{t+s}(\lambda)-W_{\infty}(\lambda)|\right]+2\pp\left[|W_{t}(\lambda)-W_{\infty}(\lambda)|\right]\to 0,\nonumber
\end{eqnarray}
as $t\to+\infty$. Thus $W_{\infty}(\lambda)$ is the $L^{1}(\pp)$-limit of $\pp\left[W_{t+s}(\lambda)|\mathcal{F}_{t}\right]$ as $t\to+\infty$.
On the other hand, by the skeleton decomposition (cf. \cite[eq.(5.9)]{CRY}), we have
\begin{equation}\label{lem2.1.1}
\pp\left[W_{t+s}(\lambda)|\mathcal{F}_{t}\right]=\e^{-c(\lambda)(t+s)}\pp\left[\langle e_{\lambda},X_{t+s}\rangle |\mathcal{F}_{t}\right]=\e^{-c(\lambda)(t+s)}\left[\langle P^{a^{*}}_{s}e_{\lambda},X_{t}\rangle+\langle P^{1}_{s}e_{\lambda},Z_{t}\rangle-\langle P^{a^{*}}_{s}e_{\lambda},Z_{t}\rangle\right],
\end{equation}
where $a^{*}=-\psi'(1)<0$. Since $P_{s}e_{\lambda}(x)=\Pi_{0}[\e^{-\lambda(B_{s}+x)}]=\e^{\lambda^{2}s/2}e_{\lambda}(x)$, we can continue the calculation in \eqref{lem2.1.1} to get that
\begin{equation*}
\pp\left[W_{t+s}(\lambda)|\mathcal{F}_{t}\right]=\e^{(-1+a^{*})s}W_{t}(\lambda)+Z_{t}(\lambda)-\e^{(-1+a^{*})s}Z_{t}(\lambda).
\end{equation*}
This implies that $\pp\left[W_{t+s}(\lambda)|\mathcal{F}_{t}\right]$ converges a.s. to $Z_{\infty}(\lambda)$ as $t\to+\infty$ and then $s\to+\infty$.  Hence we get $Z_{\infty}(\lambda)=W_{\infty}(\lambda)$ $\pp$-a.s.\qed

\begin{lemma}\label{lem2.2}
Suppose the assumptions of Lemma \ref{lem2.1} hold. Then for $|\lambda|<\sqrt{2}$, $\theta\in\R$ and $\mu\in\mf$,
\begin{equation}\label{lem2.2.1}
\pp_{\mu}\left[\e^{i\theta W_{\infty}(\lambda)}\right]=\exp\{\langle\qq_{\delta_{\cdot}}\left[\e^{i\theta Z_{\infty}(\lambda)}\right]-1,\mu\rangle\}.
\end{equation}
\end{lemma}

\proof Since for any $x\in \R$, $((X_{t})_{t\ge 0},\pp_{\delta_{x}})$ is equal in law with $((X_{t}+x)_{t\ge 0},\pp)$, one has $(W_{\infty}(\lambda),\pp_{\delta_x})\stackrel{d}{=}(e_{\lambda}(x)W_{\infty}(\lambda),\pp)$. Similarly, one has $(Z_{\infty}(\lambda),\pp_{\delta_x})\stackrel{d}{=}(e_{\lambda}(x)Z_{\infty}(\lambda),\pp)$.
Thus by Lemma \ref{lem2.1},
$$\pp_{\delta_x}\left[\e^{i\theta W_{\infty}(\lambda)}\right]=\pp\left[\e^{i\theta e_{\lambda}(x)W_{\infty}(\lambda)}\right]=\pp\left[\e^{i\theta e_{\lambda}(x)Z_{\infty}(\lambda)}\right]=\pp_{\delta_x}\left[\e^{i\theta Z_{\infty}(\lambda)}\right]\quad\forall \theta,x\in\R.$$
 By decomposing $Z_{\infty}(\lambda)$ into contributions derived from the population at time $0$, one can represent $Z_{\infty}(\lambda)$ under $\pp_{\delta_{x}}$ as
$$Z_{\infty}(\lambda)\stackrel{d}{=}\sum_{i=1}^{N_{0}}Z^{(i)}_{\infty}(\lambda)$$
where $N_{0}$ is a Poisson random variable with mean $1$, $Z^{(i)}_{\infty}(\lambda)$ are independent copies of $(Z_{\infty}(\lambda),\qq_{\delta_{x}})$ and are independent of $N_{0}$.
So we have for every $x\in\R$.
\begin{equation}\nonumber
\pp_{\delta_{x}}\left[\e^{i\theta Z_{\infty}(\lambda)}\right]=\pp_{\delta_{x}}\left[\qq_{\delta_{x}}[\e^{i\theta Z_{\infty}(\lambda)}]^{N_{0}}\right]=\exp\{\qq_{\delta_{x}}\left[\e^{i\theta Z_{\infty}(\lambda)}\right]-1\}.
\end{equation}
Thus \eqref{lem2.2.1} follows from the above identity and an application of the branching property of super-Brownian motion.\qed

\section{Proofs of main results}

\subsection{Proofs of Theorem \ref{them3.1}(i) and (ii) }
We assume (A0) and (A1) hold in this subsection. Other conditions needed will be stated explicitly.
(A1) implies that $\psi''(0+)$ exists and equals $\sigma^{2}=2b+\int_{(0,+\infty)}r^{2}\pi(dr)$.
Recall that $L$ is the integer-valued random variable with the generating function $F$ given by \eqref{eq2}.
It follows that $F''(1-)=\psi''(0+)/q<+\infty$, and thus $\mathrm{E}[L^{2}]<+\infty$.
It then follows by Lemma \ref{lem2.0} that for any $|\lambda|<1$ and $x\in\R$, $Z_{\infty}(\lambda)$ is the $L^{2}(\qq_{\delta_{x}})$-limit of $Z_{t}(\lambda)$.

\begin{lemma}
For $|\lambda|<1$ and $x\in\R$,
$$\qq_{\delta_{x}}\left[Z_{\infty}(\lambda)^{2}\right]=\frac{\sigma^{2}}{1-\lambda^{2}}e_{2\lambda}(x).$$
\end{lemma}

\proof
For any $f\in\mathcal{B}^{+}(\R)$, the second moment of $\langle f,Z_{t}\rangle$, if exists, can be expressed by
\begin{equation}\nonumber
\qq_{\delta_{x}}\left[\langle f,Z_{t}\rangle^{2}\right]=P^{1}_{t}(f^{2})(x)+\sigma^{2}\int_{0}^{t}P^{1}_{s}\left((P^{1}_{t-s}f)^{2}\right)(x)ds.
\end{equation}
It follows that
\begin{eqnarray}
\qq_{\delta_{x}}\left[Z_{t}(\lambda)^{2}\right]&=&\e^{-2c(\lambda)t}\qq_{\delta_{x}}\left[\langle e_{\lambda},Z_{t}\rangle^{2}\right]\nonumber\\
&=&\e^{-2c(\lambda)t}\left[P^{1}_{t}e_{2\lambda}(x)+\sigma^{2}\int_{0}^{t}P^{1}_{s}\left[(P^{1}_{t-s}e_{\lambda})^{2}\right](x)ds\right]\nonumber\\
&=&\e^{-(1-\lambda^{2})t}e_{2\lambda}(x)+\frac{\sigma^{2}}{1-\lambda^{2}}\left(1-\e^{-(1-\lambda^{2})t}\right)e_{2\lambda}(x).\nonumber
\end{eqnarray}
By letting $t\to+\infty$ we get $\qq_{\delta_{x}}\left[Z_{\infty}(\lambda)^{2}\right]=\frac{\sigma^{2}}{1-\lambda^{2}}e_{2\lambda}(x)$.\qed

\medskip

\noindent\textbf{Proof of Theorem \ref{them3.1}(i) and (ii):}
Define $d_{1}(t):=\e^{\frac{1-\lambda^{2}}{2}t}$ in the situation of Theorem \ref{them3.1}(i) and $d_{1}(t):=t^{1/4}\e^{t/4}$ in the situation of Theorem \ref{them3.1}(ii).
It follows by the Markov property and Lemma \ref{lem2.2} that for every $\theta\in\R$,
\begin{eqnarray}
\pp\left[\e^{i\theta d_{1}(t)(W_{\infty}(\lambda)-W_{t}(\lambda))}\right]&=&\pp\left[\e^{-i\theta d_{1}(t)W_{t}(\lambda)}\pp_{X_{t}}\left[\e^{i\theta d_{1}(t)\e^{-c(\lambda)t}W_{\infty}(\lambda)}\right]\right]\nonumber\\
&=&\pp\left[\exp\{-i\theta d_{1}(t)\e^{-c(\lambda)t}\langle e_{\lambda},X_{t}\rangle +\qq_{\delta_{\cdot}}\left[\e^{i\theta d_{1}(t)\e^{-c(\lambda)t}Z_{\infty}(\lambda)}\right]-1,X_{t}\rangle\}\right]\nonumber\\
&=&\pp\left[\exp\{\langle\qq_{\delta_{\cdot}}\left[\e^{i\theta d_{2}(t)Z_{\infty}(\lambda)}\right]-i\theta d_{2}(t)e_{\lambda}-1,X_{t} \rangle\}\right],\label{eq2.5}
\end{eqnarray}
where $d_{2}(t):=d_{1}(t)\e^{-c(\lambda)t}$.
Under our assumptions, for $|\lambda|\le \frac{\sqrt{2}}{2}$, $Z_{\infty}(\lambda)$ is the $L^{2}(\qq_{\delta_{x}})$-limit of $Z_{t}(\lambda)$ and thus $\qq_{\delta_{x}}\left[Z_{\infty}(\lambda)\right]=e_{\lambda}(x)$ and $\qq_{\delta_{x}}\left[Z_{\infty}(\lambda)^{2}\right]=\frac{\sigma^{2}}{1-\lambda^{2}}e_{2\lambda}(x)$.
One can write
\begin{equation}\label{eq2.6}
\pp\left[\e^{i\theta d_{1}(t)(W_{\infty}(\lambda)-W_{t}(\lambda))}\right]=\pp\left[\e^{\langle I_{\theta}(t,\cdot)+II_{\theta}(t,\cdot),X_{t}\rangle}\right]
\end{equation}
where
$$I_{\theta}(t,x):=\qq_{\delta_{x}}\left[\e^{i\theta d_{2}(t)Z_{\infty}(\lambda)}-1-i\theta d_{2}(t)Z_{\infty}(\lambda)+\frac{1}{2}\theta^{2}d_{2}(t)^{2}Z_{\infty}(\lambda)^{2}\right]$$
and
$$II_{\theta}(t,x):=-\frac{1}{2}\theta^{2}d_{2}(t)^{2}\qq_{\delta_{x}}\left[Z_{\infty}(\lambda)^{2}\right]=-\frac{1}{2}\theta^{2}\frac{\sigma^{2}}{1-\lambda^{2}}\,d_{2}(t)^{2}e_{2\lambda}(x).$$
Consequently we have
$$\langle II_{\theta}(t,\cdot),X_{t}\rangle=-\frac{1}{2}\theta^{2}\frac{\sigma^{2}}{1-\lambda^{2}}d_{2}(t)^{2}\e^{c(2\lambda)t}W_{t}(2\lambda).$$
In the situation of Theorem \ref{them3.1}(i), $d_{2}(t)^{2}\e^{c(2\lambda)t}\equiv 1$, and thus
\begin{equation}\label{them3.1(i)1}
\langle II_{\theta}(t,\cdot),X_{t}\rangle=-\frac{1}{2}\theta^{2}\frac{\sigma^{2}}{1-\lambda^{2}}W_{t}(2\lambda) \to -\frac{1}{2}\theta^{2}\frac{\sigma^{2}}{1-\lambda^{2}}W_{\infty}(2\lambda)\quad\pp\mbox{-a.s. as }t\to+\infty.
\end{equation}
In the situation of Theorem \ref{them3.1}(ii), $d_{2}(t)^{2}\e^{c(2\lambda)t}=\sqrt{t}$ for $\lambda=\pm\frac{\sqrt{2}}{2}$, and thus by \eqref{HRSthem1.1}
\begin{equation}\label{them3.1(ii)1}
\langle II_{\theta}(t,\cdot),X_{t}\rangle=-\theta^{2}\sigma^{2}\sqrt{t}W_{t}(\pm \sqrt{2})\to -\theta^{2}\sigma^{2}\sqrt{\frac{2}{\pi}}\partial W_{\infty}^{\pm}\mbox{ in probability as }t\to+\infty.
\end{equation}
Next we shall deal with $\langle I_{\theta}(t,\cdot),X_{t}\rangle$.
Using the fact that
$$\left|\e^{i r}-(1+i r-\frac{1}{2}r^{2})\right|\le |r|^{2}\left(\frac{|r|}{6}\wedge 1\right)\quad\forall r\in\R,$$
we have
\begin{equation}\label{eq2.9}
|I_{\theta}(t,x)|\le \theta^{2}d_{2}(t)^{2}\qq_{\delta_{x}}\left[Z^{2}_{\infty}\left(\frac{|\theta|d_{2}(t)Z_{\infty}(\lambda)}{6}\wedge 1\right)\right].
\end{equation}
In the situation of Theorem \ref{them3.1}(i), $d_{2}(t)=\e^{-\left(\lambda^{2}+\frac{1}{2}\right)t}=\e^{-c(2\lambda)t/2}$ and we have
\begin{eqnarray}
|I_{\theta}(t,x)|&\le&\theta^{2}\e^{-c(2\lambda)t}\qq_{\delta_{x}}\left[Z_{\infty}(\lambda)^{2}\left(\frac{|\theta|\e^{-c(2\lambda)t/2}Z_{\infty}(\lambda)}{6}\wedge 1\right)\right]\nonumber\\
&=&\theta^{2}\e^{-c(2\lambda)t}e_{2\lambda}(x)\qq_{\delta_{0}}\left[Z_{\infty}(\lambda)^{2}\left(\frac{|\theta|\e^{-c(2\lambda)t/2}e_{\lambda}(x)Z_{\infty}(\lambda)}{6}\wedge 1\right)\right]\nonumber\\
&=:&\theta^{2}\e^{-c(2\lambda)t}e_{2\lambda}(x)f_{\theta}(t,x).\nonumber
\end{eqnarray}
Here in the first equality we use the fact that $(Z_{\infty}(\lambda),\qq_{\delta_{x}})\stackrel{d}{=}(e_{\lambda}(x)Z_{\infty}(\lambda),\qq_{\delta_{0}})$.
Thus
\begin{equation}\label{eq2.10}
\pp\left[\left|\langle I_{\theta}(t,\cdot),X_{t}\rangle\right|\right]\le \pp\left[\langle |I_{\theta}(t,\cdot)|,X_{t}\rangle\right]\le \theta^{2}\e^{-c(2\lambda)t}\pp\left[\langle e_{2\lambda}(\cdot)f_{\theta}(t,\cdot),X_{t}\rangle\right].
\end{equation}
We have
\begin{eqnarray}
\pp\left[\langle e_{2\lambda}(\cdot)f_{\theta}(t,\cdot),X_{t}\rangle\right]&=&\e^{t}\Pi_{0}\left[e_{2\lambda}(B_{t})f_{\theta}(t,B_{t})\right]\nonumber\\
&=&\e^{c(2\lambda)t}\Pi_{0}\left[\e^{-2\lambda B_{t}-2\lambda^{2}t}f_{\theta}(t,B_{t})\right]\nonumber\\
&=&\e^{c(2\lambda)t}\Pi_{0}\left[f_{\theta}(t,B_{t}-2\lambda t)\right]\nonumber\\
&=&\e^{c(2\lambda)t}\Pi_{0}\left[f_{\theta}(t,\sqrt{t}B_{1}-2\lambda t)\right].\nonumber
\end{eqnarray}
The third equality follows from Girsanov's transform. Putting this back to \eqref{eq2.10}, we get that
\begin{equation}\label{eq2.11}
\pp\left[\left|\langle I_{\theta}(t,\cdot),X_{t}\rangle\right|\right]\le \theta^{2}\Pi_{0}\left[f_{\theta}(t,\sqrt{t}B_{1}-2\lambda t)\right].
\end{equation}
Since for every $z\in\R$ and $|\lambda|<\frac{\sqrt{2}}{2}$,
$$f_{\theta}(t,\sqrt{t}z-2\lambda t)=\qq_{\delta_{0}}\left[Z_{\infty}(\lambda)^{2}\left(\frac{|\theta|\e^{-\left(\frac{1}{2}-|\lambda|^{2}\right)t-\sqrt{t}\lambda z}Z_{\infty}(\lambda)}{6}\wedge 1\right)\right]\to 0\mbox{ as }t\to+\infty,$$
it follows by the bounded convergence theorem that $\Pi_{0}\left[f_{\theta}(t,\sqrt{t}B_{1}-2\lambda t)\right]\to 0$ as $t\to+\infty$. Thus by \eqref{eq2.11},
\begin{equation}\label{them3.1(i)2}
\langle I_{\theta}(t,\cdot),X_{t}\rangle\to 0 \mbox{ in }L^{1}(\pp)\mbox{ as }t\to+\infty.
\end{equation}

In the situation of Theorem \ref{them3.1}(ii), we may and do assume $\lambda=-\frac{\sqrt{2}}{2}$. The case $\lambda=\frac{\sqrt{2}}{2}$ can be analyzed by simple consideration of symmetry.
For notational simplification we omit $\lambda$ and write $Z_{\infty}$ for $Z_{\infty}(\lambda)$.
For $\lambda=-\frac{\sqrt{2}}{2}$, $d_{2}(t)=t^{1/4}\e^{-t}$ and \eqref{eq2.9} yields that
$$\left|I_{\theta}(t,x)\right|\le \theta^{2}t^{1/2}\e^{-2t}e_{-\sqrt{2}}(x)g_{\theta}(t,x),$$
where
$$g_{\theta}(t,x):=\qq_{\delta_{0}}\left[Z_{\infty}^{2}\left(\frac{\theta t^{1/4}\e^{-t+\frac{\sqrt{2}}{2}x}Z_{\infty}}{6}\wedge 1\right)\right].$$
Define
$$I^{1}_{\theta}(t,x):=\theta^{2}t^{1/2}\e^{-2t}e_{-\sqrt{2}}(x)g_{\theta}(t,x)1_{(-\infty,m(t)]}(x)\mbox{ and }I^{2}_{\theta}(t,x):=\theta^{2}t^{1/2}\e^{-2t}e_{-\sqrt{2}}(x)1_{(m(t),+\infty)}(x).$$
Note that for very $t\ge 0$, $x\mapsto g_{\theta}(t,x)$ is a non-decreasing function on $\R$.
It follows that
\begin{equation}\nonumber
\langle I^{1}_{\theta}(t,\cdot),X_{t}\rangle \le \theta^{2}g_{\theta}(t,m(t))t^{1/2}\e^{-2 t}\langle e_{-\sqrt{2}},X_{t}\rangle=\theta^{2}g_{\theta}(t,m(t))\sqrt{t}W_{t}(-\sqrt{2}).
\end{equation}
Then by \eqref{HRSthem1.1} and the fact that
$$g_{\theta}(t,m(t))=\qq_{\delta_{0}}\left[Z_{\infty}^{2}\left(\frac{|\theta|t^{-1/2}Z_{\infty}}{6}\wedge 1\right)\right]\to 0\mbox{ as }t\to+\infty,$$
one gets
\begin{equation}\label{them3.2.1}
\langle I^{1}_{\theta}(t,\cdot),X_{t}\rangle\to 0\mbox{ in probability as }t\to+\infty.
\end{equation}
On the other hand,
\begin{eqnarray}
\langle I^{2}_{\theta}(t,\cdot),X_{t}\rangle&=&\theta^{2}t^{1/2}\e^{-2t}\langle e_{-\sqrt{2}}1_{(m(t),+\infty)},X_{t}\rangle\nonumber\\
&=&\theta^{2}t^{1/2}\e^{-2 t}\langle e_{-\sqrt{2}}(\cdot+m(t))1_{(m(t),+\infty)}(\cdot+m(t)),X_{t}-m(t)\rangle\nonumber\\
&=&\theta^{2}t^{-1}\langle e_{-\sqrt{2}}1_{(0,+\infty)},\mathcal{E}_{t}\rangle.\label{eq3.9}
\end{eqnarray}
We note that (A3) and \eqref{condi4} hold under our assumptions of Theorem \ref{them3.1}(ii). Then one can easily verify that the function $\e_{-\sqrt{2}}(x)1_{(0,+\infty)}(x)$ satisfies the assumptions of Lemma \ref{lem4.1}. Hence by Lemma \ref{lem4.1} $\langle e_{-\sqrt{2}}1_{(0,+\infty)},\mathcal{E}_{t}\rangle$ converges in distribution to a finite random variable. Thus by \eqref{eq3.9}
$\langle I^{2}_{\theta}(t,\cdot),X_{t}\rangle $ converges in distribution (or, equivalently, converges in probability) to $0$. This together with \eqref{them3.2.1}  and the fact that
$|\langle I_{\theta}(t,\cdot),X_{t}\rangle|\le \langle I^{1}_{\theta}(t,\cdot),X_{t}\rangle+\langle I^{2}_{\theta}(t,\cdot),X_{t}\rangle$ yields that
\begin{equation}\label{them3.1(ii)2}
\langle I_{\theta}(t,\cdot),X_{t}\rangle \to 0\mbox{ in probability as }t\to+\infty.
\end{equation}
Therefore by \eqref{them3.1(i)1} and \eqref{them3.1(i)2}, we have
$$\langle I_{\theta}(t,\cdot)+II_{\theta}(t,\cdot),X_{t}\rangle\to -\frac{1}{2}\theta^{2}\,\frac{\sigma^{2}}{1-\lambda^{2}}\,W_{\infty}(2\lambda)\mbox{ in probability as }t\to+\infty,$$
in the situation of Theorem \ref{them3.1}(i), and by \eqref{them3.1(ii)1} and \eqref{them3.1(ii)2},
$$\langle I_{\theta}(t,\cdot)+II_{\theta}(t,\cdot),X_{t}\rangle\to -\theta^{2}\sigma^{2}\sqrt{\frac{2}{\pi}}\,\partial W_{\infty}^{\pm}\mbox{ in probability as }t\to+\infty$$
 in the situation of Theorem \ref{them3.1}(ii).
In either situation, one can show by \eqref{eq2.6} that the characteristic function of $d_{1}(t)(W_{\infty}(\lambda)-W_{t}(\lambda))$ converges pointwisely to the characteristic function of the claimed limit. Hence we complete the proof.\qed

\subsection{Proofs of Theorem \ref{them5.1}(i) and (ii)}

We assume (A0) and (A2) hold in this section. Other conditions needed will be stated explicitly.
Recall that $\bar{\pi}(r)=\int_{(r,+\infty)}\pi(dr)$ for $r>0$ and $\psi_{0}(\lambda)=\lambda+\psi(\lambda)$ for $\lambda\ge 0$. By Fubini's theorem one may write $\psi_{0}(\lambda)$ as
\begin{eqnarray}
\psi_{0}(\lambda)&=&b\lambda^{2}+\lambda\int_{0}^{+\infty}\left(1-\e^{-\lambda r}\right)\bar{\pi}(r)dr\nonumber\\
&=&\lambda\left(b\lambda+\int_{0}^{+\infty}\left(1-\e^{-\lambda r}\right)\bar{\pi}(r)dr\right)\nonumber\\
&=:&\lambda g(\lambda).
\end{eqnarray}
Note that $g(\lambda)$ is a Bernstein function, and thus \cite[Proposition 3.6]{SSV} allows us to extend $\psi_{0}$ (and then $\psi$) analytically onto the right complex plane $\C^{+}_{0}:=\{z\in\C:\ \re z>0\}$, and continuously onto the closure $\C^{+}:=\{z\in\C:\ \re z\ge 0\}$. We still use $\psi_{0}$ (and $\psi$) to denote such an extension.

Recall that $A_{1},\beta$ and $\delta$ are the constants given in (A2). Define another measure on $(0,+\infty)$ by
$$\pi_{1}(dr):=A_{1}(1+\beta)r^{-(2+\beta)}1_{\{r>0\}}dr,$$
and correspondingly $\bar{\pi}_{1}(r):=\int_{(r,+\infty)}\pi_{1}(dy)=A_{1}r^{-1-\beta}$ for $r>0$.
Then we define for $z\in \C^{+}$,
\begin{eqnarray}
\psi_{1}(z)&:=&z\int_{0}^{+\infty}\left(1-\e^{-z r}\right)\bar{\pi}_{1}(r)dr\nonumber\\
&=&A_{2}z^{1+\beta},\nonumber
\end{eqnarray}
where $A_{2}:=A_{1}\Gamma(1-\beta)/\beta$. Here in the last equality we use the fact that $z^{\beta}=\frac{\beta}{\Gamma(1-\beta)}\int_{0}^{+\infty}\left(1-\e^{-z r}\right)r^{-(1+\beta)}dr$ for $z\in\C^{+}$. Obviously the function $\psi_{1}$ is holomorphic on $\C^{+}_{0}$ and continuous on $\C^{+}$.

\begin{lemma}\label{lem5.1}
For every $\epsilon\in (0,\delta]$, there are constants $c_{1},c_{2}>0$ such that
$$|\psi_{0}(z)-\psi_{1}(z)|\le c_{1}|z|^{2}+c_{2}|z|^{1+\beta+\epsilon}\quad\forall z\in \C^{+}.$$
\end{lemma}

\proof Fix an arbitrary $\epsilon\in (0,\delta]$. For $z\in \C^{+}$,
\begin{eqnarray}
|\psi_{0}(z)-\psi_{1}(z)|&=&\left|bz^{2}+z\int_{0}^{+\infty}\left(1-\e^{-z r}\right)(\bar{\pi}(r)-\bar{\pi}_{1}(r))dr\right|\nonumber\\
&\le&b|z|^{2}+|z|\int_{0}^{+\infty}(|z|r\wedge 2)|\bar{\pi}(r)-\bar{\pi}_{1}(r)|dr.\label{lem5.1.1}
\end{eqnarray}
Condition (A2) yields that there are constants $M\ge 1$ and $c_{3}>0$ such that
$$|\bar{\pi}(r)-\bar{\pi}_{1}(r)|\le \frac{c_{3}}{r^{1+\beta+\delta}}\quad\forall r\ge M.$$
This implies in particular that for $\epsilon\in (0,\delta]$,
$$|\bar{\pi}(r)-\bar{\pi}_{1}(r)|\le \frac{c_{3}}{r^{1+\beta+\epsilon}}\quad\forall r\ge M.$$
Now we write
$$\int_{0}^{+\infty}(|z|r\wedge 2)|\bar{\pi}(r)-\bar{\pi}_{1}(r)|dr=I(z)+II(z),$$
where
$$I(z):=\int_{0}^{M}(|z|r\wedge 2)|\bar{\pi}(r)-\bar{\pi}_{1}(r)|dr\mbox{ and }II(z):=\int_{M}^{+\infty}(|z|r\wedge 2)|\bar{\pi}(r)-\bar{\pi}_{1}(r)|dr.$$
We have
\begin{eqnarray}
II(z)&\le&c_{3}\int_{M}^{+\infty}(|z|r\wedge 2)\frac{1}{r^{1+\beta+\epsilon}}\,dr\nonumber\\
&=&c_{3}\left(1_{\{2/|z|>M\}}\int_{M}^{2/|z|}\frac{|z|r}{r^{1+\beta+\epsilon}}dr+\int_{2/|z|\vee M}^{+\infty}\frac{2}{r^{1+\beta+\epsilon}}dr\right)\nonumber\\
&=&c_{3}\left(1_{\{2/|z|>M\}}(1-\beta-\epsilon)|z|\left(\left(\frac{2}{|z|}\right)^{1-\beta-\epsilon}-M^{1-\beta-\epsilon}\right)+\frac{2}{\beta+\epsilon}\left(\frac{2}{|z|}\vee M\right)^{-\beta-\epsilon}\right)\nonumber\\
&\le&c_{4}|z|^{\beta+\epsilon}.\label{lem5.1.2}
\end{eqnarray}
On the other hand, by Fubini's theorem
\begin{eqnarray}
I(z)&\le&\int_{0}^{M}|z|r(\bar{\pi}(r)+\bar{\pi}_{1}(r))dr\nonumber\\
&=&|z|\left(\frac{1}{2}\int_{(0,+\infty)}\left(r^{2}\wedge M^{2}\right)\pi(dr)+\frac{1}{2}\int_{(0,+\infty)}\left(r^{2}\wedge M^{2}\right)\pi_{1}(dr)\right)\nonumber\\
&=:&c_{5}|z|.\label{lem5.1.3}
\end{eqnarray}
Putting \eqref{lem5.1.2} and \eqref{lem5.1.3} back into \eqref{lem5.1.1}, we prove the assertion.\qed

The above lemma implies in particular that for every $\eta\in (0,\beta]$, there is some constant $c_{1}>0$ such that
\begin{equation}\label{cortolem5.1}
|\psi_{0}(z)|\le c_{1}\left(|z|^{2}+|z|^{1+\eta}\right)\quad\forall z\in \C^{+}.
\end{equation}
In fact, Lemma \ref{lem5.1} yields that for $\epsilon\in (0,\delta]$,
\begin{equation}\nonumber
|\psi_{0}(z)|\lesssim |z|^{2}+|z|^{1+\beta+\epsilon}+|z|^{1+\beta}\asymp |z|^{2}\vee |z|^{1+\beta+\epsilon}\vee|z|^{1+\beta}=|z|^{2}\vee|z|^{1+\beta}\quad\forall z\in \C^{+}.
\end{equation}
Here for nonnegative functions $f$ and $g$ we write $f(x)\lesssim g(x)$ if there exists a positive constant $c_{1}$ such that $f(x)\le c_{1}g(x)$ for all $x$, and write $f(x)\asymp g(x)$ if $f(x)\lesssim g(x)$ and $g(x)\lesssim f(x)$.
Then we get from the above inequality that for $\eta\in (0,\beta]$,
$$|\psi_{0}(z)|\lesssim |z|^{2}\vee|z|^{1+\eta}\asymp |z|^{2}+|z|^{1+\eta}\quad\forall z\in\C^{+}.$$

The generating function $F$ of $L$ can be extended analytically onto the disk $B_{1}:=\{z\in\C: |z|<1\}$ and continuously onto the closure $\bar{B}_{1}:=\{z\in\C:\ |z|\le 1\}$ by the representation:
$$F(z)=\mathrm{E}[z^{L}]\quad\forall z\in \bar{B}_{1}.$$
It necessarily holds that
\begin{equation}\nonumber
q(F(z)-z)=\psi(1-z)\quad\forall z\in \bar{B}_{1}.
\end{equation}
Here $\psi(z)$ denotes the extension of $\psi(\lambda)$ onto $\C^{+}$ which is homomorphic on $\C^{+}_{0}$ and continuous on $\C^{+}$.

For $\lambda\in\R$, define
$$W^{\lambda}(\theta,x):=\qq_{\delta_{x}}\left[e^{i\theta Z_{\infty}(\lambda)}\right]-1\quad\forall \theta,x\in\R.$$
We notice that $Z_{t}(0)=\e^{-t}\|Z_{t}\|$ where $((\|Z_{t}\|)_{t\ge 0},\qq_{\delta_{x}})$ is a continuous-time Galton-Watson process whose distributions are independent of the starting point $x$. Thus the distribution of $Z_{\infty}(0)$ under $\qq_{\delta_{x}}$ is independent of $x$, and $W^{0}(\theta,x)=W^{0}(\theta,0)$ for all $x\in \R$. We simply denote $W^{0}(\theta,x)$ by $W^{0}(\theta)$.

\begin{lemma}\label{lem5.4}
For $|\lambda|<\sqrt{2}$,
\begin{equation}\label{lem5.4.0}
W^{\lambda}(\theta,x)=i\theta e_{\lambda}(x)+\Pi_{x}\left[\int_{0}^{+\infty}\e^{s}\psi_{0}\left(-W^{\lambda}(\e^{-c(\lambda)s}\theta,B_{s})\right)ds\right]\quad\forall \theta,x\in\R.
\end{equation}
In particular, for $\lambda=0$,
$$W^{0}(\theta)=i\theta+\int_{0}^{+\infty}\e^{s}\psi_{0}\left(-W^{0}(\e^{-s}\theta)\right)ds\quad\forall \theta\in\R.$$
\end{lemma}

\proof For $\lambda\in \R$ and $\theta\in\R$, define
$$U^{\lambda}_{\theta}(t,x):=\qq_{\delta_{x}}\left[\e^{i\theta Z_{t}(\lambda)}\right]\mbox{ and }V^{\lambda}_{\theta}(t,x):=U^{\lambda}_{\theta}(t,x)-1\quad\forall t\ge 0,\ x\in\R.$$
First we consider the case $\lambda\not=0$.
We have
\begin{eqnarray}
U^{\lambda}_{\theta}(t,x)&=&\qq_{\delta_{x}}\left[\exp\{i\theta\e^{-c(\lambda)t}\langle e_{\lambda},Z_{t}\rangle\}\right]\nonumber\\
&=&\e^{-qt}\Pi_{x}\left[\exp\{i\theta\e^{-c(\lambda)t}e_{\lambda}(B_{t})\right]+\Pi_{x}\left[\int_{0}^{t}q\e^{-q s}\sum_{k}p_{k}\qq_{\delta_{B_{s}}}\left[\exp\{i\theta\e^{-c(\lambda)s}Z_{t-s}(\lambda)\}\right]^{k}ds\right]\nonumber\\
&=&\e^{-qt}\Pi_{x}\left[\exp\{i\theta\e^{-c(\lambda)t}e_{\lambda}(B_{t})\right]+\Pi_{x}\left[\int_{0}^{t}q\e^{-q s}F\left(U^{\lambda}_{\theta\e^{-c(\lambda)s}}(t-s,B_{s})\right)ds\right].\label{lem5.4.1}
\end{eqnarray}
Here in the second equality, the first term corresponds to the case when the ancestor of the branching Brownian motion is still alive at time $t$, and the second term corresponds to the case when it dies at $s\in (0,t]$.
Since $(Z_{r}(\lambda),\qq_{\delta_{y+\frac{c(\lambda)}{\lambda}s}})\stackrel{d}{=}(\e^{-c(\lambda)s}Z_{r}(\lambda),\qq_{\delta_{y}})$, we have
$$U^{\lambda}_{\theta\e^{-c(\lambda)s}}(r,y)=U^{\lambda}_{\theta}(r,y+\frac{c(\lambda)}{\lambda}s)\quad\forall r\ge0,\ y\in\R.$$
It then follows by \eqref{lem5.4.1} that
\begin{eqnarray}
U^{\lambda}_{\theta}(t,x)&=&\e^{-qt}\Pi_{x}\left[\exp\{i\theta e_{\lambda}(B_{t}+\frac{c(\lambda)}{\lambda}t)\right]+\Pi_{x}\left[\int_{0}^{t}q\e^{-q s}F\left(U^{\lambda}_{\theta}(t-s,B_{s}+\frac{c(\lambda)}{\lambda}s)\right)ds\right]\nonumber\\
&=&\e^{-qt}\Pi_{x}\left[\exp\{i\theta e_{\lambda}(\xi^{\lambda}_{t})\right]+\Pi_{x}\left[\int_{0}^{t}q\e^{-q s}F\left(U^{\lambda}_{\theta}(t-s,\xi^{\lambda}_{s})\right)ds\right],\label{eq2.24}
\end{eqnarray}
where $\xi^{\lambda}_{t}:=B_{t}+c(\lambda)t/\lambda$ for $t\ge 0$.
Immediately, we get
$$\re U^{\lambda}_{\theta}(t,x)=\e^{-qt}\Pi_{x}\left[\re\exp\{i\theta e_{\lambda}(\xi^{\lambda}_{t})\right]+\Pi_{x}\left[\int_{0}^{t}q\e^{-q s}\re F\left(U^{\lambda}_{\theta}(t-s,\xi^{\lambda}_{s})\right)ds\right],$$
and
$$\im U^{\lambda}_{\theta}(t,x)=\e^{-qt}\Pi_{x}\left[\im\exp\{i\theta e_{\lambda}(\xi^{\lambda}_{t})\right]+\Pi_{x}\left[\int_{0}^{t}q\e^{-q s}\im F\left(U^{\lambda}_{\theta}(t-s,\xi^{\lambda}_{s})\right)ds\right].$$
Applying the integral identity established in \cite[Lemma A.1]{EKW} to the above identities, one gets that
$$\re U^{\lambda}_{\theta}(t,x)=\Pi_{x}\left[\re\exp\{i\theta e_{\lambda}(\xi^{\lambda}_{t})\right]+\Pi_{x}\left[\int_{0}^{t}q\left(\re F\left(U^{\lambda}_{\theta}(t-s,\xi^{\lambda}_{s})\right)-\re U^{\lambda}_{\theta}(t-s,\xi^{\lambda}_{s})\right)ds\right],$$
and
$$\im U^{\lambda}_{\theta}(t,x)=\Pi_{x}\left[\im\exp\{i\theta e_{\lambda}(\xi^{\lambda}_{t})\right]+\Pi_{x}\left[\int_{0}^{t}q\left(\im F\left(U^{\lambda}_{\theta}(t-s,\xi^{\lambda}_{s})\right)-\im U^{\lambda}_{\theta}(t-s,\xi^{\lambda}_{s})\right)ds\right].$$
Combining the above two identities, we have
\begin{eqnarray}
 U^{\lambda}_{\theta}(t,x)&=&\Pi_{x}\left[\exp\{i\theta e_{\lambda}(\xi^{\lambda}_{t})\right]+\Pi_{x}\left[\int_{0}^{t}q\left(F\left(U^{\lambda}_{\theta}(t-s,\xi^{\lambda}_{s})\right)- U^{\lambda}_{\theta}(t-s,\xi^{\lambda}_{s})\right)ds\right]\nonumber\\
 &=&\Pi_{x}\left[\exp\{i\theta e_{\lambda}(\xi^{\lambda}_{t})\right]+\Pi_{x}\left[\int_{0}^{t}\psi\left(1-\left(U^{\lambda}_{\theta}(t-s,\xi^{\lambda}_{s})\right)\right)ds\right],\nonumber
\end{eqnarray}
or equivalently,
$$V^{\lambda}_{\theta}(t,x)=\Pi_{x}\left[\exp\{i\theta e_{\lambda}(\xi^{\lambda}_{t})\right]-1+\Pi_{x}\left[\int_{0}^{t}\psi_{0}\left(-V^{\lambda}_{\theta}(t-s,\xi^{\lambda}_{s})\right)+V^{\lambda}_{\theta}(t-s,\xi^{\lambda}_{s})ds\right].$$
Again, by applying \cite[Lemma A.1]{EKW} to the real and imaginary parts of the above equation and combining the identities together, one can show that
\begin{eqnarray}
V^{\lambda}_{\theta}(t,x)
%&=&\e^{t}\left(\Pi_{x}\left[\exp\{i\theta e_{\lambda}(\xi^{\lambda}_{t})\right]-1\right)+\Pi_{x}\left[\int_{0}^{t}\e^{s}\psi_{0}\left(-V^{\lambda}_{\theta}(t-s,\xi^{\lambda}_{s})\right)ds\right]\nonumber\\
&=&I(t,x)+i\theta e_{\lambda}(x)+\Pi_{x}\left[\int_{0}^{t}\e^{s}\psi_{0}\left(-V^{\lambda}_{\theta}(t-s,\xi^{\lambda}_{s})\right)ds\right],\label{eq2.25}
\end{eqnarray}
where $I(t,x):=\e^{t}\left(\Pi_{x}\left[\exp\{i\theta e_{\lambda}(\xi^{\lambda}_{t})\right]-1\right)-i\theta e_{\lambda}(x)$.
Note that $W^{\lambda}(\theta,x)=\lim_{t\to+\infty}V^{\lambda}_{\theta}(t,x)$ and the fact that $W^{\lambda}(\e^{-c(\lambda)s}\theta,y)=W^{\lambda}(\theta,y+\frac{c(\lambda)}{\lambda}s)$.  We can prove \eqref{lem5.4.0} by letting $t\to+\infty$ in \eqref{eq2.25} once we prove that
\begin{equation}\label{lem5.4.2}
\lim_{t\to+\infty}I(t,x)=0\quad\forall x\in \R.
\end{equation}
We note that for any $1\le p\le 2$,
\begin{equation}\nonumber
\e^{t}\Pi_{x}\left[e_{\lambda}(\xi^{\lambda}_{t})^{p}\right]=\e^{t-c(\lambda)pt}\Pi_{x}[e_{\lambda}(B_{t})]=\e^{-c_{p}(\lambda)t}e_{\lambda p}(x)
\end{equation}
where $c_{p}(\lambda):=c(\lambda)p-c(\lambda p)=(p-1)(1-\lambda^{2}p/2)$. In particular, $\e^{t}\Pi_{x}\left[e_{\lambda}(\xi^{\lambda}_{t})\right]=e_{\lambda}(x)$.
For $|\lambda|<\sqrt{2}$ we can take $p\in (1,2/\lambda^{2})\cap (1,2]\not=\emptyset$ such that $c_{p}(\lambda)>0$. It follows by Lemma \ref{lemA.2} that
\begin{eqnarray}
|I(t,x)|&=&\e^{t}\left|\Pi_{x}\left[\e^{i\theta e_{\lambda}(\xi^{\lambda}_{t})}-1-i\theta e_{\lambda}(\xi^{\lambda}_{t})\right]\right|\nonumber\\
&\le&c_{1}|\theta|^{p}\e^{t}\Pi_{x}[e_{\lambda}(\xi^{\lambda}_{t})^{p}]\nonumber\\
&=&c_{1}|\theta|^{p}\e^{-c_{p}(\lambda)t}e_{\lambda p}(x)\to 0\mbox{ as }t\to+\infty.
\end{eqnarray}
Thus we prove \eqref{lem5.4.2}.

The proof for the case $\lambda=0$ is similar. Since the distribution of $\|Z_{t}\|$ under $\qq_{\delta_{x}}$ is independent of the starting point $x$, we define for all $x\in\R$,
$$U_{\theta}(t):=\qq_{\delta_{x}}[\e^{i\theta \|Z_{t}\|}]\mbox{ and }V_{\theta}(t):=U_{\theta}(t)-1,\quad\forall t\ge 0,\theta\in\R.$$
Similarly as in \eqref{lem5.4.1}, we have
$$U_{\theta}(t)=\e^{-qt}\e^{i\theta}+\int_{0}^{t}q\e^{-q s}F(U_{\theta}(t-s))ds.$$
 Again, by applying the integral identities of \cite[Lemma A.1]{EKW} to the real and imaginary parts of the above equation, and then to the real and imaginary parts of the equation for $V_{\theta}(t)$, one can show that
 $$V_{\theta}(t)=\e^{t}(e^{i\theta}-1)+\int_{0}^{t}\e^{s}\psi_{0}(-V_{\theta}(t-s))ds.$$
Consequently, we have
$$V_{\theta\e^{-t}}(t)=\e^{t}(e^{i\theta\e^{-t}}-1)+\int_{0}^{t}\e^{s}\psi_{0}(-V_{\theta\e^{-t}}(t-s))ds.$$
Note that $W^{0}(\theta)=\lim_{t\to +\infty}V_{\theta\e^{-t}}(t)$. By letting $t\to+\infty$ in the above equation we prove the second assertion.\qed

\begin{lemma}\label{lem5.2}
For $1<p<1+\beta$, $\int_{(1,+\infty)}r^{p}\pi(dr)<+\infty$.
\end{lemma}
\proof By Fubini's theorem we have
\begin{eqnarray}
\int_{(1,+\infty)}r^{p}\pi(dr)&=&\int_{(1,+\infty)}\pi(dr)\int_{0}^{r}py^{p-1}dy\nonumber\\
&=&p\int_{0}^{+\infty}y^{p-1}\bar{\pi}(y\vee 1)dy\nonumber\\
&=&p\left(\bar{\pi}(1)\int_{0}^{1}y^{p-1}dy+\int_{1}^{+\infty}y^{p-1}\bar{\pi}(y)dy\right).
\end{eqnarray}
We only need to show the final integral of the above equation is finite. Condition (A2) yields that there are constants $c_{1},c_{2}>0$ such that for $y$ sufficiently large,
$$y^{p-1}\bar{\pi}(y)\le y^{p-1}(c_{1}y^{-1-\beta}+c_{2}y^{-1-\beta-\delta})=c_{1}y^{-1-(1+\beta-p)}+c_{2}y^{-1-(1+\beta+\delta-p)}.$$
Thus $\int_{1}^{+\infty}y^{p-1}\bar{\pi}(y)dy$ is finite for $1<p<1+\beta$. \qed

\medskip

\noindent\textbf{Proof of Theorem \ref{them5.1}(i) and (ii):} Let $d_{3}(t):=\e^{\frac{\beta}{2}(\frac{2}{1+\beta}-\lambda^{2})t}$ in the situation of Theorem \ref{them5.1}(i), and
$d_{3}(t):=t^{\frac{1}{2(1+\beta)}}\e^{\frac{\beta^{2}}{(1+\beta)^{2}}t}$ in the situation of Theorem \ref{them5.1}(ii).
We have showed in \eqref{eq2.5} that
\begin{equation}\label{eq3.25}
\pp\left[\e^{i\theta d_{3}(t)(W_{\infty}(\lambda)-W_{t}(\lambda))}\right]
=\pp\left[\exp\{\langle\qq_{\delta_{\cdot}}\left[\e^{i\theta d_{4}(t)Z_{\infty}(\lambda)}\right]-i\theta d_{4}(t)e_{\lambda}(\cdot)-1,X_{t}\rangle\}\right],\quad\forall \theta\in\R
\end{equation}
where $d_{4}(t):=d_{3}(t)\e^{-c(\lambda)t}$.
It follows by Lemma \ref{lem5.4} that
\begin{eqnarray*}
\qq_{\delta_{x}}[\e^{i\theta d_{4}(t)Z_{\infty}(\lambda)}]-i\theta d_{4}(t)\e_{\lambda}(x)-1
&=&W^{\lambda}(\theta d_{4}(t),x)-i\theta d_{4}(t)e_{\lambda}(x)\\
&=&\Pi_{x}\left[\int_{0}^{+\infty}\e^{s}\psi_{0}\left(-W^{\lambda}(\theta d_{4}(t)\e^{-c(\lambda)s},B_{s})\right)ds\right]\\
&=&I_{\theta}(t,x)+II_{\theta}(t,x)+III_{\theta}(t,x),
\end{eqnarray*}
where
\begin{eqnarray*}
I_{\theta}(t,x)&:=&\Pi_{x}\left[\int_{0}^{+\infty}\e^{s}\psi_{1}\left(-i\theta d_{4}(t)\e^{-c(\lambda)s}e_{\lambda}(B_{s})\right)ds\right],\\
II_{\theta}(t,x)&:=&\Pi_{x}\left[\int_{0}^{+\infty}\e^{s}\left(\psi_{1}\left(-W^{\lambda}(\theta d_{4}(t)\e^{-c(\lambda)s},B_{s})\right)-\psi_{1}\left(-i\theta d_{4}(t)\e^{-c(\lambda)s}e_{\lambda}(B_{s})\right)\right)ds\right],\\
III_{\theta}(t,x)&:=&\Pi_{x}\left[\int_{0}^{+\infty}\e^{s}\left(\psi_{0}\left(-W^{\lambda}(\theta d_{4}(t)\e^{-c(\lambda)s},B_{s})\right)-\psi_{1}\left(-W^{\lambda}(d_{4}(t)\e^{-c(\lambda)s},B_{s})\right)\right)ds\right],\\
\end{eqnarray*}
In the remaining of this proof, we shall show that
\begin{description}
\item{(i)} In the situation of Theorem \ref{them5.1}(i), for $|\lambda|<\frac{\sqrt{2}}{1+\beta}$,
$$\langle I_{\theta}(t,\cdot),X_{t}\rangle\to c_{\lambda,\beta}(-i\theta)^{1+\beta}W_{\infty}((1+\beta)\lambda)\quad\pp\mbox{-a.s. as }t\to+\infty,$$
and in the situation of Theorem \ref{them5.1}(ii), for $|\lambda|=\pm\frac{\sqrt{2}}{1+\beta}$,
$$\langle I_{\theta}(t,\cdot),X_{t}\rangle\to c_{\beta}(-i\theta)^{1+\beta}\sqrt{\frac{2}{\pi}}\partial W_{\infty}^{\pm}\mbox{ in probability as }t\to+\infty.$$
\item{(ii)} $\langle II_{\theta}(t,\cdot),X_{t}\rangle\to 0$ in probability as $t\to+\infty$.
\item{(iii)} $\langle III_{\theta}(t,\cdot),X_{t}\rangle\to 0$ in probability as $t\to+\infty$.
\end{description}
 Using the above convergence results and \eqref{eq3.25}, one can easily show that the characteristic function of $d_{3}(t)(W_{\infty}(\lambda)-W_{t}(\lambda))$ converges pointwisely to that of the claimed limit as $t\to+\infty$ and hence complete the proof.

\medskip

\noindent Proof of (i): Recall that $\psi_{1}(z)=A_{2}z^{1+\beta}$ for $z\in\C^{+}$ where $A_{2}=A_{1}\Gamma(1-\beta)/\beta$. We have
\begin{eqnarray}
I_{\theta}(t,x)
&=&A_{2}d_{4}(t)^{1+\beta}(-i\theta)^{1+\beta}\int_{0}^{+\infty}\e^{s-c(\lambda)(1+\beta)s}\Pi_{x}\left[e_{\lambda}(B_{s})^{1+\beta}\right]ds\nonumber\\
&=&A_{2}d_{4}(t)^{1+\beta}(-i\theta)^{1+\beta}\int_{0}^{+\infty}\e^{-c_{1+\beta}(\lambda)s}e_{(1+\beta)\lambda}(x)ds.\nonumber
\end{eqnarray}
Here $c_{1+\beta}(\lambda)=c(\lambda)(1+\beta)-c(\lambda(1+\beta))=\beta(1-\frac{\lambda^{2}}{2}(1+\beta))$. For $|\lambda|\le \frac{\sqrt{2}}{1+\beta}$, $c_{1+\beta}(\lambda)>0$.  Hence we can continue the above calculation to get that
$$I_{\theta}(t,x)=c_{\lambda,\beta}(-i\theta)^{1+\beta}d_{4}(t)^{1+\beta}e_{(1+\beta)\lambda}(x)$$
where $c_{\lambda,\beta}=A_{2}c_{1+\beta}(\lambda)=A_{2}\beta^{-1}(1-\frac{1}{2}\lambda^{2}(1+\beta))^{-1}$.

In the situation of Theorem \ref{them5.1}(i), for $|\lambda|<\frac{\sqrt{2}}{1+\beta}$, $d_{4}(t)=\exp\{-\left(\frac{1+\beta}{2}\lambda^{2}+\frac{1}{1+\beta}\right)t\}=\exp\{-\frac{c((1+\beta)\lambda)}{1+\beta}t\}$, and thus
\begin{eqnarray}
\langle I_{\theta}(t,\cdot),X_{t}\rangle&=&c_{\lambda,\beta}(-i\theta)^{1+\beta}\e^{-c(\lambda(1+\beta))t}\langle e_{(1+\beta)\lambda},X_{t}\rangle\nonumber\\
&=&c_{\lambda,\beta}(-i\theta)^{1+\beta}W_{t}((1+\beta)\lambda)\nonumber\\
&\to&c_{\lambda,\beta}(-i\theta)^{1+\beta}W_{\infty}((1+\beta)\lambda)\quad\pp\mbox{-a.s. as }t\to+\infty.\nonumber
\end{eqnarray}
On the other hand, in the situation of Theorem \ref{them5.1}(ii), 
\eqref{condi3} is satisfied by Lemma \ref{lem5.2}, and thus $\sqrt{t}W_{t}(\pm\sqrt{2})\to \sqrt{2/\pi}\partial W^{\pm}_{\infty}$ in probability.
For $\lambda=\pm\frac{\sqrt{2}}{1+\beta}$,
$d_{4}(t)=t^{\frac{1}{2(1+\beta)}}\e^{-\frac{2}{1+\beta}t}$, and thus
\begin{eqnarray}
\langle I_{\theta}(t,\cdot),X_{t}\rangle&=&c_{\lambda,\beta}(-i\theta)^{1+\beta}\sqrt{t}\e^{-2t}\langle e_{\pm\sqrt{2}},X_{t}\rangle\nonumber\\
&=&c_{\lambda,\beta}(-i\theta)^{1+\beta}\sqrt{t}W_{t}(\pm \sqrt{2})\nonumber\\
&\to& c_{\beta}(-i\theta)^{1+\beta}\partial W_{\infty}^{\pm}\mbox{ in probability as }t\to +\infty,\nonumber
\end{eqnarray}
where $c_{\beta}=c_{\lambda,\beta}\sqrt{2/\pi}=A_{1}\beta^{-3}(1+\beta)\Gamma(1-\beta)\sqrt{2/\pi}$.

\medskip

\noindent Proof of (ii): We have
\begin{eqnarray}
|II_{\theta}(t,x)|&\le&\Pi_{x}\left[\int_{0}^{+\infty}\e^{s}\left|\psi_{1}\left(-W^{\lambda}(\theta d_{4}(t)\e^{-c(\lambda)s},B_{s})\right)-\psi_{1}\left(-i\theta d_{4}(t)\e^{-c(\lambda)s}e_{\lambda}(B_{s})\right)\right|ds\right]\nonumber\\
&=&A_{2}\Pi_{x}\left[\int_{0}^{+\infty}\e^{s}\left|\left(-W^{\lambda}(\theta d_{4}(t)\e^{-c(\lambda)s},B_{s})\right)^{1+\beta}-\left(-i\theta d_{4}(t)\e^{-c(\lambda)s}e_{\lambda}(B_{s})\right)^{1+\beta}\right|ds\right].\label{them5.1.1}
\end{eqnarray}
We notice that $\re W^{\lambda}(\theta,y)=\qq_{\delta_{y}}\left[\cos(\theta Z_{\infty}(\lambda))\right]-1\le 0$ for all $\theta,y\in\R$, and thus $-W^{\lambda}(\theta,y)\in \C^{+}$. Using the fact that (cf. \cite[Lemma A.3]{RSSZ})
$$\left|z_{1}^{1+\beta}-z_{2}^{1+\beta}\right|\le (1+\beta)\left(|z_{1}|^{\beta}+|z_{2}|^{\beta}\right)\left|z_{1}-z_{2}\right|\quad\forall z_{1},z_{2}\in \C^{+},$$
we have for all $\theta,y\in\R$
\begin{equation}
\left|(-W^{\lambda}(\theta,y))^{1+\beta}-(-i\theta e_{\lambda}(y))^{1+\beta}\right|\le(1+\beta)\left(|W^{\lambda}(\theta,y)|^{\beta}+|i\theta e_{\lambda}(y)|^{\beta}\right)\left|W^{\lambda}(\theta,y)-i\theta e_{\lambda}(y)\right|.\label{them5.1.2}
\end{equation}
By Lemma \ref{lem5.4},
\begin{eqnarray}
\left|W^{\lambda}(\theta,y)-i\theta e_{\lambda}(y)\right|&=&\left|\Pi_{y}\left[\int_{0}^{+\infty}\e^{s}\psi_{0}\left(-W^{\lambda}(\theta\e^{-c(\lambda)s},B_{s})\right)ds\right]\right|\nonumber\\
&\le&\Pi_{y}\left[\int_{0}^{+\infty}\e^{s}\left|\psi_{0}\left(-W^{\lambda}(\theta\e^{-c(\lambda)s},B_{s})\right)\right|ds\right].\label{them5.1.4}
\end{eqnarray}
Let $\epsilon_{1}\in (0,\beta)$. It follows by \eqref{cortolem5.1} that there is some constant $c_{1}>0$ such that $|\psi_{0}(z)|\le c_{1}(|z|^{2}+|z|^{1+\epsilon_{1}})$ for all $z\in \C^{+}$.
Hence by \eqref{them5.1.4}
\begin{equation}
|W^{\lambda}(\theta,y)-i\theta e_{\lambda}(y)|\le c_{1}\int_{0}^{+\infty}\e^{s}\left(\Pi_{y}\left[|W^{\lambda}(\theta \e^{-c(\lambda)s},B_{s})|^{2}\right]+\Pi_{y}\left[|W^{\lambda}(\theta \e^{-c(\lambda)s},B_{s})|^{1+\epsilon_{1}}\right]\right)ds.\label{eq4.17}
\end{equation}
Note that for $\theta,y\in\R$,
\begin{equation}\label{them5.1.5}
|W^{\lambda}(\theta,y)|=\left|\qq_{\delta_{y}}\left[\e^{i\theta Z_{\infty}(\lambda)}\right]-1\right|\le \qq_{\delta_{y}}\left[\left|\e^{i\theta Z_{\infty}(\lambda)}-1\right|\right]\le\qq_{\delta_{y}}\left[2\wedge |\theta|Z_{\infty}(\lambda)\right]\le 2\wedge |\theta|e_{\lambda}(y).
\end{equation}
Consequently we have $\Pi_{y}\left[\left|W^{\lambda}(\theta\e^{-c(\lambda)s},B_{s})\right|^{1+\epsilon_{1}}\right]
\le\Pi_{y}\left[|\theta|^{1+\epsilon_{1}}\e^{-c(\lambda)s(1+\epsilon_{1})}e_{\lambda}(B_{s})^{1+\epsilon_{1}}\right]$. Moreover, we have
\begin{eqnarray}
\Pi_{y}\left[|W^{\lambda}(\theta\e^{-c(\lambda)s},B_{s})|^{2}\right]&\le&\Pi_{y}\left[4\wedge |\theta|^{2}\e^{-2c(\lambda)s}e_{\lambda}(B_{s})^{2}\right]\nonumber\\
&=&\Pi_{y}\left[|\theta|^{2}\e^{-2c(\lambda)s}e_{\lambda}(B_{s})^{2};|\theta|\e^{-c(\lambda)s}e_{\lambda}(B_{s})\le 2\right]+4\Pi_{y}\left[|\theta|\e^{-c(\lambda)s}e_{\lambda}(B_{s})>2\right]\nonumber\\
&\le&2^{1-\epsilon_{1}}\Pi_{y}\left[|\theta|^{1+\epsilon_{1}}\e^{-c(\lambda)(1+\epsilon_{1})s}e_{\lambda}(B_{s})^{1+\epsilon_{1}};|\theta|\e^{-c(\lambda)s}e_{\lambda}(B_{s})\le 2\right]
\nonumber\\
&&\quad+4\Pi_{y}\left[\frac{|\theta|^{1+\epsilon_{1}}\e^{-c(\lambda)(1+\epsilon_{1})s}e_{\lambda}(B_{s})^{1+\epsilon_{1}}}{2^{1+\epsilon_{1}}};|\theta|\e^{-c(\lambda)s}e_{\lambda}(B_{s})>2\right]\nonumber\\
&\le&2^{1-\epsilon_{1}}\Pi_{y}\left[|\theta|^{1+\epsilon_{1}}\e^{-c(\lambda)(1+\epsilon_{1})s}e_{\lambda}(B_{s})^{1+\epsilon_{1}}\right].\label{eq4.18}
\end{eqnarray}
Then we get from \eqref{eq4.17} that
\begin{eqnarray}
\left|W^{\lambda}(\theta,y)-i\theta e_{\lambda}(y)\right|&\le&c_{2}|\theta|^{1+\epsilon_{1}}\int_{0}^{+\infty}\e^{s-c(\lambda)(1+\epsilon_{1})s}\Pi_{y}\left[e_{\lambda}(B_{s})^{1+\epsilon_{1}}\right]ds\nonumber\\
&=&c_{2}|\theta|^{1+\epsilon_{1}}\int_{0}^{+\infty}\e^{-c_{1+\epsilon_{1}}(\lambda)s}e_{(1+\epsilon_{1})\lambda}(y)ds\nonumber\\
&=&c_{3}|\theta|^{1+\epsilon_{1}}e_{(1+\epsilon_{1})\lambda}(y).\nonumber
\end{eqnarray}
In the last equality we use the fact that $c_{1+\epsilon_{1}}(\lambda)=\epsilon_{1}\left(1-\frac{\lambda^{2}(1+\epsilon_{1})}{2}\right)>0$ for $1+\epsilon_{1}\in (1,1+\beta)\subset (1,2/\lambda^{2})$.
It follows by the above inequality, \eqref{them5.1.5} and \eqref{them5.1.2} that
$$\left|(-W^{\lambda}(\theta,y))^{1+\beta}-(-i\theta e_{\lambda}(y))^{1+\beta}\right|\le c_{4}|\theta|^{1+\beta+\epsilon_{1}}e_{(1+\beta+\epsilon_{1})\lambda}(y)\quad\forall \theta,y\in\R.$$
Putting this back into \eqref{them5.1.1} we get that
\begin{eqnarray}
|II_{\theta}(t,x)|&\le&c_{4}A_{2}\Pi_{x}\left[\int_{0}^{+\infty}\e^{s}|\theta d_{4}(t)\e^{-c(\lambda)s}|^{1+\beta+\epsilon_{1}}e_{(1+\beta+\epsilon_{1})\lambda}(B_{s})ds\right]\nonumber\\
&=&c_{5}|\theta|^{1+\beta+\epsilon_{1}}d_{4}(t)^{1+\beta+\epsilon_{1}}\int_{0}^{+\infty}\e^{-c_{1+\beta+\epsilon_{1}}(\lambda)s}e_{(1+\beta+\epsilon_{1})\lambda}(x)ds.\nonumber
\end{eqnarray}
For $|\lambda|\le \frac{\sqrt{2}}{1+\beta}$, one may take $\epsilon_{1}>0$ so small that $1+\beta+\epsilon_{1}<\frac{2}{\lambda^{2}}$, and thus $c_{1+\beta+\epsilon_{1}}(\lambda)=(\beta+\epsilon_{1})\left(1-\frac{\lambda^{2}(1+\beta+\epsilon_{1})}{2}\right)>0$. Then we get from the above inequalities that
\begin{equation}\label{them5.1.6}
|II_{\theta}(t,x)|\le c_{6}|\theta|^{1+\beta+\epsilon_{1}}d_{4}(t)^{1+\beta+\epsilon_{1}}e_{(1+\beta+\epsilon_{1})\lambda}(x).
\end{equation}
In the situation of Theorem \ref{them5.1}(i), the above inequality yields that
$$\left|\langle II_{\theta}(t,\cdot),X_{t}\rangle\right|\le c_{6}|\theta|^{1+\beta+\epsilon_{1}}\e^{-l_{\lambda,\beta}(\epsilon_{1})t}W_{t}((1+\beta+\epsilon_{1})\lambda),$$
where $l_{\lambda,\beta}(\epsilon_{1}):=c((1+\beta)\lambda)\left(1+\frac{\epsilon_{1}}{1+\beta}\right)=\epsilon_{1}\left(\frac{1}{1+\beta}
-\frac{1+\beta+\epsilon_{1}}{2}\lambda^{2}\right)$.
For $|\lambda|<\frac{\sqrt{2}}{1+\beta}$,  $\frac{1}{1+\beta}-\frac{1+\beta}{2}\lambda^{2}>0$, and one may take $\epsilon_{1}>0$ so small that $l_{\lambda,\beta}(\epsilon_{1})>0$. Thus we get from the above inequality that $\langle II_{\theta}(t,\cdot),X_{t}\rangle\to 0\quad \pp$-a.s. as $t\to+\infty$.

In the situation of Theorem \ref{them5.1}(ii), without loss of generality we assume $\lambda=-\frac{\sqrt{2}}{1+\beta}$.
By \eqref{them5.1.6} we have
$$|II_{\theta}(t,x)|\le c_{6}|\theta|^{1+\beta+\epsilon_{1}}\left(t^{-1/2}\e^{-2t}\right)^{1+\frac{\epsilon_{1}}{1+\beta}}e_{-\sqrt{2}\left(1+\frac{\epsilon_{1}}{1+\beta}\right)}(x).$$
Thus
\begin{eqnarray}
|\langle II_{\theta}(t,\cdot),X_{t}\rangle|&\le&c_{6}|\theta|^{1+\beta+\epsilon_{1}}\left(t^{-1/2}\e^{-2t}\right)^{1+\frac{\epsilon_{1}}{1+\beta}}
\langle e_{-\sqrt{2}\left(1+\frac{\epsilon_{1}}{1+\beta}\right)},X_{t}\rangle\nonumber\\
&=&c_{6}|\theta|^{1+\beta+\epsilon_{1}}\left(t^{-1/2}\e^{-2t}\right)^{1+\frac{\epsilon_{1}}{1+\beta}}
\langle e_{-\sqrt{2}\left(1+\frac{\epsilon_{1}}{1+\beta}\right)}(\cdot+m(t)),X_{t}-m(t)\rangle\nonumber\\
&=&c_{6}|\theta|^{1+\beta+\epsilon_{1}}t^{-\left(1+\frac{\epsilon_{1}}{1+\beta}\right)}\langle e_{-\sqrt{2}\left(1+\frac{\epsilon_{1}}{1+\beta}\right)},\mathcal{E}_{t}\rangle.\nonumber
\end{eqnarray}
In view of Lemma \ref{lem5.2}, (A3) and \eqref{condi4} hold under our assumptions of Theorem \ref{them5.1}(ii). Then one can easily verify that the function $e_{-\sqrt{2}\left(1+\frac{\epsilon_{1}}{1+\beta}\right)}(x)$ satisfies the conditions of Lemma \ref{lem4.1}.
It follows that $\langle e_{-\sqrt{2}\left(1+\frac{\epsilon_{1}}{1+\beta}\right)},\mathcal{E}_{t}\rangle$ converges in distribution to some finite random variable as $t\to+\infty$. Hence we get by the above inequality that $\langle II_{\theta}(t,\cdot),X_{t}\rangle$ converges in distribution (and thus converges in probability) to $0$ as $t\to+\infty$.

\medskip

\noindent Proof of (iii): Let $\epsilon_{2}\in (0,\delta)$. We have by Lemma \ref{lem5.1},
\begin{eqnarray}
|III_{\theta}(t,x)|&\le&\Pi_{x}\left[\int_{0}^{+\infty}\e^{s}\left|\psi_{0}\left(-W^{\lambda}(\theta d_{4}(t)\e^{-c(\lambda)s},B_{s})\right)-\psi_{1}\left(-W^{\lambda}(\theta d_{4}(t)\e^{-c(\lambda)s},B_{s})\right)\right|ds\right]\nonumber\\
&\le&c_{1}\Pi_{x}\left[\int_{0}^{+\infty}\e^{s}\left(\left|W^{\lambda}(\theta d_{4}(t)\e^{-c(\lambda)s},B_{s})\right|^{2}+\left|W^{\lambda}(\theta d_{4}(t)\e^{-c(\lambda)s},B_{s})\right|^{1+\beta+\epsilon_{2}}\right)ds\right].\label{them5.1.3}
\end{eqnarray}
By \eqref{them5.1.5}, $\Pi_{x}\left[\left|W^{\lambda}(\theta \e^{-c(\lambda)s},B_{s})\right|^{1+\beta+\epsilon_{2}}\right]\le \Pi_{x}\left[|\theta|^{1+\beta+\epsilon_{2}}\e^{-c(\lambda)(1+\beta+\epsilon_{2})s}e_{\lambda}(B_{s})^{1+\beta+\epsilon_{2}}\right]$ for every $\theta\in\R$ and $s\ge 0$.
Similarly as in \eqref{eq4.18}, one can show that
$$\Pi_{x}\left[|W^{\lambda}(\theta\e^{-c(\lambda)s},B_{s})|^{2}\right]\le 2^{1-(\beta+\epsilon_{2})}\Pi_{x}\left[|\theta|^{1+\beta+\epsilon_{2}}\e^{-c(\lambda)(1+\beta+\epsilon_{2})s}e_{\lambda}(B_{s})^{1+\beta+\epsilon_{2}}\right]\quad\forall \theta\in\R,\ s\ge 0.$$
Thus we get from \eqref{them5.1.3} that
\begin{eqnarray}
|III_{\theta}(t,x)|&\le&c_{2}|\theta|^{1+\beta+\epsilon_{2}}d_{4}(t)^{1+\beta+\epsilon_{2}}
\int_{0}^{+\infty}\e^{s}\Pi_{x}\left[\e^{-c(\lambda)(1+\beta+\epsilon_{2})s}e_{\lambda}(B_{s})^{1+\beta+\epsilon_{2}}\right]\nonumber\\
&=&c_{2}|\theta|^{1+\beta+\epsilon_{2}}d_{4}(t)^{1+\beta+\epsilon_{2}}\int_{0}^{+\infty}\e^{-c_{1+\beta+\epsilon_{2}}(\lambda)s}e_{(1+\beta+\epsilon_{2})\lambda}(x)ds\nonumber\\
&=&c_{3}|\theta|^{1+\beta+\epsilon_{2}}d_{4}(t)^{1+\beta+\epsilon_{2}}e_{(1+\beta+\epsilon_{2})\lambda}(x).\label{eq3.33}
\end{eqnarray}
Here to get the final equality we take $\epsilon_{2}$ so small that $c_{1+\beta+\epsilon_{2}}(\lambda)>0$.
Applying similar argument as in the end of the proof of (ii) (with \eqref{them5.1.6} replaced by \eqref{eq3.33}),  one can show that
$\langle III_{\theta}(t,\cdot),X_{t}\rangle\to 0$ $\pp$-a.s. in the situation of Theorem \ref{them5.1}(i) and in probability in the situation of Theorem \ref{them5.1}(ii).
Therefore we complete the proof.\qed

\subsection{Proofs of Theorem \ref{them3.1}(iii) and Theorem \ref{them5.1}(iii) }

\begin{lemma}\label{lem4.2}
Assume (A3) and \eqref{condi4} hold. For $\frac{\sqrt{2}}{2}<|\lambda|<\sqrt{2}$, if there exists some $p\in (\frac{\sqrt{2}}{|\lambda|},2]$ such that $\qq_{\delta_{0}}\left[Z_{\infty}(\lambda)^{p}\right]<+\infty$, then there is a finite and non-degenerate random variable $\eta_{\lambda}$ such that
\begin{equation}\label{lem4.2.0}
t^{\frac{3|\lambda|}{2\sqrt{2}}}\e^{\frac{1}{2}\left(\sqrt{2}-|\lambda|\right)^{2}t}\left(W_{\infty}(\lambda)-W_{t}(\lambda)\right)\stackrel{d}{\to} \eta_{\lambda}\mbox{ as }t\to+\infty.
\end{equation}
Moreover, the characteristic function of $\eta_{\lambda}$ is given by \eqref{eq:cf}
\end{lemma}

\proof We give proof for $\lambda<0$. The case $\lambda>0$ can be analyzed similarly by symmetry. Let
\begin{equation*}
d_{5}(t):=t^{\frac{3|\lambda|}{2\sqrt{2}}}\e^{\frac{1}{2}\left(\sqrt{2}-|\lambda|\right)^{2}t}=\e^{\left(\frac{1}{2}\lambda^{2}+1\right)t+\lambda\left(\sqrt{2} t-\frac{3}{2\sqrt{2}}\log t\right)}=\e^{c(\lambda)t+\lambda m(t)}.
\end{equation*}
For $\theta\in\R$, we have
\begin{eqnarray}
\pp\left[\e^{i\theta d_{5}(t)(W_{\infty}(\lambda)-W_{t}(\lambda))}\right]&=&\pp\left[\exp\{-i\theta d_{5}(t)W_{t}(\lambda)+\langle\qq_{\delta_{\cdot}}\left[\e^{i\theta d_{5}(t)\e^{-c(\lambda)t}Z_{\infty}(\lambda)}\right]-1,X_{t}\rangle\}\right]\nonumber\\
&=&\pp\left[\exp\{-i\theta \e^{\lambda m(t)}\}\langle e_{\lambda},X_{t}\rangle+\langle\qq_{\delta_{\cdot}}\left[\e^{i\theta \e^{\lambda m(t)}Z_{\infty}(\lambda)}\right]-1,X_{t}\rangle\right]\nonumber\\
&=&\pp\left[\exp\{-i\theta\langle e_{\lambda}(\cdot-m(t)),X_{t}\rangle+\langle\qq_{\delta_{\cdot-m(t)}}\left[\e^{i\theta Z_{\infty}(\lambda)}\right]-1,X_{t}\rangle\}\right]\nonumber\\
&=&\pp\left[\exp\{\langle\qq_{\delta_{\cdot}}\left[\e^{i\theta Z_{\infty}(\lambda)}\right]-1-i\theta e_{\lambda}(\cdot),X_{t}-m(t)\rangle\}\right]\nonumber\\
&=&\pp\left[\e^{\langle F_{\lambda,\theta},\mathcal{E}_{t}\rangle}\right],\nonumber
\end{eqnarray}
where $F_{\lambda,\theta}(x):=\qq_{\delta_{x}}\left[\e^{i\theta Z_{\infty}(\lambda)}\right]-1-i\theta e_{\lambda}(x)$ for $x\in\R$.
Here the first equality follows from the Markov property and Lemma \ref{lem2.2}.
We shall show that
$\pp\left[\e^{\langle F_{\lambda,\theta},\mathcal{E}_{t}\rangle}\right]$ converges pointwisely to the characteristic function of some non-degenerate random variable.
Under our assumptions, for $|\lambda|<\sqrt{2}$, $W_{\infty}(\lambda)$ is the $L^{1}(\pp_{\delta_{x}})$-limit of $W_{t}(\lambda)$ for every $x\in\R$, and thus $\pp_{\delta_{x}}[W_{\infty}(\lambda)]=e_{\lambda}(x)$. It then follows by Lemma \ref{lem2.2} that
$$\qq_{\delta_{x}}\left[Z_{\infty}(\lambda)\right]=\pp_{\delta_{x}}\left[W_{\infty}(\lambda)\right]=e_{\lambda}(x).$$
Thus we have $F_{\lambda,\theta}(x)=\qq_{\delta_{x}}\left[\e^{i\theta Z_{\infty}(\lambda)}-1-i\theta Z_{\infty}(\lambda)\right]$.
By Lemma \ref{lemA.2}, for $\frac{\sqrt{2}}{2}<|\lambda|<\sqrt{2}$ and $p\in (\frac{\sqrt{2}}{|\lambda|},2]$ such that $\qq_{\delta_{0}}\left[Z_{\infty}(\lambda)^{p}\right]<+\infty$,
\begin{equation*}
|F_{\lambda,\theta}(x)|\le c_{1}\qq_{\delta_{x}}\left[|\theta|^{p}Z_{\infty}(\lambda)^{p}\right]=c_{1}|\theta|^{p}\e^{-\lambda p x}\qq_{\delta_{0}}\left[Z_{\infty}(\lambda)^{p}\right]=c_{2}|\theta|^{p}\e^{|\lambda|p x}\quad\forall x\in\R,
\end{equation*}
where $c_{1}=2^{3-2p}$ and $c_{2}=c_{1}\qq_{\delta_{0}}\left[Z_{\infty}(\lambda)^{p}\right]$.
Since $|\lambda|p>\sqrt{2}$, one can easily verify by the above inequality that $|F_{\lambda,\theta}(x)|$ is a locally bounded function satisfying $|F_{\lambda,\theta}(x)|1_{(-\infty,0]}(x)\in\mathcal{H}$. Consequently, the assumptions of Lemma \ref{lem4.1} are satisfied by the functions $\re^{\pm}F_{\lambda,\theta}(x)$, $\im^{\pm}F_{\lambda,\theta}(x)$ and their linear combinations. Hence by Lemma \ref{lem4.1} and the Cram\'{e}r-Wold device (cf. \cite[Corollary 5.5]{Kallenberg}),
\begin{multline}
\left(\langle \re^{+}F_{\lambda,\theta},\mathcal{E}_{t}\rangle,\langle \re^{-}F_{\lambda,\theta},\mathcal{E}_{t}\rangle,\langle \im^{+}F_{\lambda,\theta},\mathcal{E}_{t}\rangle,\langle \im^{-}F_{\lambda,\theta},\mathcal{E}_{t}\rangle\right)\\
\stackrel{d}{\to}\left(\langle \re^{+}F_{\lambda,\theta},\mathcal{E}_{\infty}\rangle,\langle \re^{-}F_{\lambda,\theta},\mathcal{E}_{\infty}\rangle,\langle \im^{+}F_{\lambda,\theta},\mathcal{E}_{\infty}\rangle,\langle \im^{-}F_{\lambda,\theta},\mathcal{E}_{\infty}\rangle\right)
\end{multline}
as $t\to+\infty$, where the limit is a vector of finite random variables. This yields that for every $\theta\in\R$,
$$\pp\left[\e^{i\theta d_{5}(t)(W_{\infty}(\lambda)-W_{t}(\lambda))}\right]=\pp\left[\e^{\langle F_{\lambda,\theta},\mathcal{E}_{t}\rangle}\right]\to \mathrm{E}\left[\e^{\langle F_{\lambda,\theta},\mathcal{E}_{\infty}\rangle}\right]=:g_{\lambda}(\theta)\mbox{ as }t\to+\infty.$$
Here we remark that $\re F_{\lambda,\theta}(x)=\qq_{\delta_{x}}\left[\cos(\theta Z_{\infty}(\lambda))\right]-1\le 0$ for all $x\in\R$, and thus $|g_{\lambda}(\theta)|\le \mathrm{E}\left[\left|\e^{\langle F_{\lambda,\theta},\mathcal{E}_{\infty}\rangle}\right|\right]=\mathrm{E}\left[\e^{\langle \re F_{\lambda,\theta},\mathcal{E}_{\infty}\rangle}\right]<+\infty$ for all $\theta\in\R$.
Since $F_{\lambda,0}(x)\equiv 0$ by definition, $g_{\lambda}(0)=1$. This implies that $d_{5}(t)(W_{\infty}(\lambda)-W_{t}(\lambda))$ converges in distribution to some finite random variable $\eta_{\lambda}$ with characteristic function $\mathrm{E}\left[\e^{i\theta\eta_{\lambda}}\right]=g_{\lambda}(\theta)$.
To show $\eta_{\lambda}$ is non-degenerate, it suffices to show that $\{\theta\in\R:\ |g_{\lambda}(\theta)|<1\}$ is a non-empty set. We take $\theta=1$. Note that by definition
$$\re F_{\lambda,1}(x)=\qq_{\delta_{x}}\left[\cos Z_{\infty}(\lambda)\right]-1=\qq_{\delta_{0}}\left[\cos(\e^{-\lambda x}Z_{\infty}(\lambda))\right]-1.$$
It is easy to see that $\re F_{\lambda,1}(x)$ is a non-positive continuous function on $\R$ with $\{x\in\R:\ \re F_{\lambda,1}(x)<0\}$ being a non-empty open set.
We know by \cite[Theorem 2.10]{RYZ} that $\langle f,\mathcal{E}_{\infty}\rangle$ vanishes if and only if the non-zero points of $f$ has zero Lebsgue measure.
Hence $\langle \re F_{\lambda,1},\mathcal{E}_{\infty}\rangle$ is negative with positive probability.  It follows that
$$|g_{\lambda}(1)|=\left|\mathrm{E}\left[\e^{\langle F_{\lambda,1},\mathcal{E}_{\infty}\rangle\rangle}\right]\right|\le\mathrm{E}\left[\left|\e^{\langle F_{\lambda,1},\mathcal{E}_{\infty}\rangle}\right|\right]=
\mathrm{E}\left[\e^{\langle \re F_{\lambda,1},\mathcal{E}_{\infty}\rangle}\right]<1.$$
By Lemma \ref{lem2.2}, for $x,\theta\in\R$,
\begin{eqnarray*}
F_{\lambda,\theta}(x)&=&\log \pp_{\delta_{x}}\left[\e^{i\theta W_{\infty}(\lambda)}\right]-i\theta e_{\lambda}(x)\\
&=&\log \left(\e^{-i\theta e_{\lambda}(x)}\pp_{\delta_{x}}\left[\e^{i\theta W_{\infty}(\lambda)}\right]\right)\\
&=&\log \pp_{\delta_{x}}\left[\e^{i\theta(W_{\infty}(\lambda)-W_{0}(\lambda))}\right].
\end{eqnarray*}
Thus the characteristic function of $\eta_{\lambda}$ can also be written as
\begin{equation}\label{eq:cf}
g_{\lambda}(\theta)=\mathrm{E}[\e^{i\theta\eta_{\lambda}}]=\mathrm{E}\left[\e^{\langle \log \pp_{\delta_{\cdot}}\left[e^{i\theta (W_{\infty}(\lambda)-W_{0}(\lambda))}\right],\mathcal{E}_{\infty}\rangle}\right].
\end{equation}
We complete the proof.\qed

\begin{proposition}\label{prop2}
Assume (A3) holds. For $\frac{\sqrt{2}}{2}<|\lambda|<\sqrt{2}$, if there exists some $p\in (\frac{\sqrt{2}}{|\lambda|},\frac{2}{\lambda^{2}})$ and $p\le 2$ such that $\int_{(1,+\infty)} r^{p}\pi(dr)<+\infty$, then \eqref{lem4.2.0} holds.
\end{proposition}
\proof Recall that $L$ has the law of the offspring distribution of the skeleton BBM. By \cite[Lemma A.1]{RYZ}, $\mathrm{E}[L^{p}]<+\infty$ for the constant $p$ satisfying our assumptions, and thus by Lemma \ref{lem2.0} $Z_{\infty}(\lambda)$ is the $L^{p}(\qq_{\delta_{0}})$-limit of $Z_{\infty}(\lambda)$. Hence the proposition follows directly from Lemma \ref{lem4.2}.\qed

\medskip

\noindent\textbf{Proof of Theorem \ref{them3.1}(iii) and Theorem \ref{them5.1}(iii):}
In the situation of Theorem \ref{them3.1}(iii), one may take $p=2$ for $\frac{\sqrt{2}}{2}<|\lambda|<1$ and $p\in (\frac{\sqrt{2}}{|\lambda|},\frac{2}{\lambda^{2}})\subset (1,2)$ for $1\le |\lambda|<\sqrt{2}$, such that $\int_{(1,+\infty)}r^{p}\pi(dr)<+\infty$.
In the situation of Theorem \ref{them5.1}(iii), for $\frac{\sqrt{2}}{1+\beta}<|\lambda|<\sqrt{2}$, by Lemma \ref{lem5.2} one may take $p\in (\frac{\sqrt{2}}{|\lambda|},2]\cap (1,(1+\beta)\wedge \frac{2}{\lambda^{2}})\not=\emptyset$ such that $\int_{(1,+\infty)}r^{p}\pi(dr)<+\infty$. Thus the assertion follows by Proposition \ref{prop2}.\qed

\appendix
\section{Appendix}

\begin{lemma}\label{lemA.0}
Let $\psi$ be given in \eqref{2.1} with $\psi(+\infty)=+\infty$.
The following statements are equivalent.
\begin{description}
\item{(i)} There are constants $\gamma\in (0,1]$, $M>0$ and $c_{1}>0$ such that
$$\psi_{0}(\lambda)\ge c_{1}\lambda^{1+\gamma}\quad\forall \lambda\ge M.$$
\item{(ii)} There are constants $\gamma\in (0,1]$, $c>0$ and $d>0$ such that
$$\psi(\lambda)\ge -c\lambda+d\lambda^{1+\gamma}\quad\forall \lambda>0.$$
\item{(iii)} Either $b>0$ or $b=0$ and there are constants $\gamma\in (0,1]$, $\epsilon>0$ and $c_{2}>0$ such that
$$\int_{0}^{1}\bar{\pi}(r)(r\wedge\lambda)dr\ge c_{2}\lambda^{1-\gamma}\quad\forall 0<\lambda<\epsilon.$$
\end{description}
\end{lemma}

\proof
Without loss of generality, we assume $\psi'(0+)=1$ and $\psi(1)=0$.

(i)$\Rightarrow$(ii): Assume (i) holds. Then immediately one has
$$\psi(\lambda)=-\lambda+\psi_{0}(\lambda)\ge -\lambda+c_{1}\lambda^{1+\gamma}\quad\forall \lambda\ge M.$$
Since $\psi(1)=0$ and $\psi'(1)>0$, we have
$$\lim_{\lambda\to 1+}\frac{-\lambda+\lambda^{1+\gamma}}{\psi(\lambda)}=\lim_{\lambda\to 1+}\frac{-1+(1+\gamma)\lambda^{\gamma}}{\psi'(\lambda)}=-\frac{\gamma}{\psi'(1)}<+\infty.$$
This implies that the function $x\mapsto (-x+x^{1+\gamma})/\psi(x)$ is bounded from above on $(1,M)$. Thus there is some constant $c_{3}>0$ such that
$$\psi(\lambda)\ge -c_{3}\lambda+c_{3}\lambda^{1+\gamma}\quad\forall 1<\lambda<M.$$
Finally for $0<\lambda\le 1$, we have
$$\psi(\lambda)\ge -\lambda\ge -2\lambda+\lambda^{1+\gamma}.$$
Combining the above results, we prove that (ii) holds for $c=c_{3}\vee 2$ and $d=c_{1}\wedge c_{3}\wedge 1$.

(ii)$\Rightarrow$(i): Assume (ii) holds. Then
$$\psi_{0}(\lambda)=\psi(\lambda)+\lambda\ge (1-c)\lambda+d\lambda^{1+\gamma}=\lambda^{1+\gamma}\left(d+\frac{1-c}{\lambda^{\gamma}}\right)\quad\forall \lambda>0.$$
Since $\lim_{\lambda\to+\infty}(1-c)/\lambda^{\gamma}=0$, it follows that for $\lambda>0$ sufficiently large,
$\psi_{0}(\lambda)\ge \frac{d}{2}\lambda^{1+\gamma}.$
Hence we prove (i).

(i)$\Leftrightarrow$(iii): Suppose $b=0$. In this case,
$$\psi_{0}(\lambda)=\int_{(0,+\infty)}(\e^{-\lambda r}-1+\lambda r)\pi(dr)=\lambda\int_{0}^{+\infty}(1-\e^{-\lambda r})\bar{\pi}(r)dr.$$
It then follows by \cite[Lemma 3.4]{SSV} that
$$\psi_{0}(\lambda)\asymp\lambda^{2}\int_{0}^{1/\lambda}\bar{\bar{\pi}}(r)dr,$$
where $\bar{\bar{\pi}}(r):=\int_{r}^{+\infty}\bar{\pi}(z)dz$.
So (i) holds if and only if there is some constants $c_{4}>0$ and $M_{1}>1$,
$$\int_{0}^{1/\lambda}\bar{\bar{\pi}}(r)dr\ge c_{4}\lambda^{-(1-\gamma)}\quad\forall \lambda\ge M_{1},$$
or equivalently,
\begin{equation}\label{lemA.0.1}
\int_{0}^{\lambda}\bar{\bar{\pi}}(r)dr\ge c_{4}\lambda^{1-\gamma}\quad\forall 0<\lambda<M_{1}^{-1}.
\end{equation}
We note that by Fubini's theorem
$$\int_{0}^{\lambda}\bar{\bar{\pi}}(r)dr=\int_{0}^{\lambda}\int_{r}^{+\infty}\bar{\pi}(z)dz dr=\int_{0}^{+\infty}\bar{\pi}(z)(z\wedge \lambda)dz
=(\int_{0}^{1}+\int_{1}^{+\infty})\bar{\pi}(z)(z\wedge\lambda)dz,$$
and that for $0<\lambda<M_{1}^{-1}$,
$$\int_{1}^{+\infty}\bar{\pi}(z)(z\wedge \lambda)dz=\lambda\int_{1}^{+\infty}\bar{\pi}(z)dz=\lambda\int_{(1,+\infty)}(z-1)\pi(dz),$$
where the final integral is finite. Therefore \eqref{lemA.0.1} holds if and only if
$$\int_{0}^{1}\bar{\pi}(z)(z\wedge\lambda)dz\ge c_{4}\lambda^{1-\gamma}\quad\forall 0<\lambda<M_{1}^{-1}.$$
Hence we complete the proof.\qed

 \begin{lemma}\label{lemA.2}
 Suppose $Y$ is a real valued random variable. For $p\in [1,2]$ such that $\mathrm{E}[|Y|^{p}]<+\infty$,
 $$\left|\mathrm{E}[\e^{iY}-iY-1]\right|\le 2^{3-2p}\mathrm{E}[|Y|^{p}].$$
 \end{lemma}
\proof Using the fact that $|e^{ix}-ix-1|\le 2|x|\wedge \frac{|x|^{2}}{2}$, we have
\begin{eqnarray*}
\left|\mathrm{E}[\e^{iY}-iY-1]\right|&\le&\mathrm{E}[\left|\e^{iY}-iY-1\right|]\\
&\le&\mathrm{E}\left[2|Y|\wedge \frac{Y^{2}}{2}\right]\\
&=&\frac{1}{2}\mathrm{E}[|Y|^{2};|Y|\le 4]+2\mathrm{E}[|Y|;|Y|>4]\\
&\le&\frac{1}{2}\mathrm{E}[|Y|^{p}\cdot|Y|^{2-p};|Y|\le 4]+2\mathrm{E}[|Y|\cdot\frac{|Y|^{p-1}}{4^{p-1}};|Y|>4]\\
&\le&2^{3-2p}\mathrm{E}[|Y|^{p};|Y|\le 4]+2^{3-2p}\mathrm{E}[|Y|^{p};|Y>4|]\\
&=&2^{3-2p}\mathrm{E}[|Y|^{p}].
\end{eqnarray*}\qed

\medskip

{ Ting Yang}

School of Mathematics and Statistics, Beijing Institute of Technology, Beijing, 100081, P.R.China;

Beijing Key Laboratory on MCAACI,

Beijing, 100081, P.R. China.

Email: yangt@bit.edu.cn

\end{document}